\documentclass[11pt,a4paper]{article}
\usepackage[utf8]{inputenc}
\usepackage{amsmath}
\usepackage{amsfonts}  
\usepackage{amsmath}
\usepackage{amssymb}
\usepackage{amsthm}
\usepackage{setspace}
\usepackage{enumerate}
\usepackage{hyperref}
\usepackage{xcolor}
\usepackage[color, leftbars]{changebar}
\usepackage{lipsum}
\usepackage[scr=rsfso]{mathalfa}
\usepackage[Symbol]{upgreek}
\usepackage{scalerel}
\usepackage{graphicx}
\usepackage{bbm}

\usepackage{fancyhdr}
\usepackage{etoolbox}
\usepackage[T1]{fontenc}

\usepackage[a4paper, right=3.2cm, left=3.2cm, top=4.0cm, bottom=3.9cm]{geometry}

\newtheoremstyle{noenddot}
{\topsep}
{\topsep}
{\itshape}
{-4pt}
{\bfseries}
{}
{ } 
{\thmname{#1}\thmnumber{ #2}\thmnote{ \normalfont(#3)}}

\theoremstyle{noenddot}
\newtheorem{theorem}{}[section]

\newcommand{\thm}{\textnormal{\textbf{Theorem.~}}}
\newcommand{\proposition}{\textnormal{\textbf{Proposition.~}}}

\newcommand{\lemma}{\textnormal{\textbf{Lemma.~}}}

\newcommand{\varupalpha}{\scalebox{0.9}[1]{\scaleobj{0.95}{\upalpha}}\hspace{0.05em}}
\newcommand{\varupbeta}{\scalebox{0.9}[1]{\scaleobj{0.95}{\upbeta}}}

\newcommand{\varupmu}{\scalebox{0.9}[1]{\scaleobj{0.95}{\upmu}}}
\newcommand{\varuppsi}{\scalebox{0.9}[1]{\scaleobj{0.92}{\uppsi}}}
\newcommand{\varupvarphi}{\scalebox{0.9}[1]{\scaleobj{0.95}{\upvarphi}}}
\newcommand{\varupphi}{\scalebox{0.9}[1]{\scaleobj{0.95}{\upphi}}}
\newcommand{\varuptheta}{\scalebox{0.9}[1]{\scaleobj{0.95}{\uptheta}}}
\newcommand{\varupnu}{\scalebox{0.9}[1]{\scaleobj{0.95}{\upnu}}}
\newcommand{\varuppi}{\scalebox{0.9}[1]{\scaleobj{0.92}{\uppi}}}

\newcommand{\T}{T}

\title{\vspace{-5ex}\large\textbf{MIXING AND EQUIPARTITION FOR AUTOMORPHISM INVARIANT PROCESSES ON REGULAR TREES}}
\author{\normalsize\textsc{Felix Pogorzelski and Elias Zimmermann}}
\date{\vspace{-5ex}}

\begin{document}
\begin{spacing}{1.1}

\maketitle

    \renewcommand{\abstractname}{\vspace{-\baselineskip}}

    \begin{abstract}
        \textsc{Abstract}. The paper is devoted to  equipartition of measured information for finite state processes over regular trees
         whose laws are invariant under all parity preserving tree automorphisms.  We show almost everywhere equipartition for ergodic processes along spheres and balls in every horosphere. Moreover, under a quantitive mixing condition we obtain a Shannon-McMillan-Breiman theorem along metric spheres of even radius. 
    \end{abstract}

\section{Introduction}

\noindent \textbf{Background.} Let $(\Omega,\mathscr{A},\varupmu)$ be a standard probability space, $E$ be a finite set and $\varupalpha = \{\varupalpha^{i}\colon i \in E\}$ be a finite partition of $\Omega$. We define the \textit{information function} corresponding to $\varupalpha$ by
    \[I(\varupalpha) := -\sum_{i \in E} \mathbbm{1}_{\upalpha_{i}}\log\varupmu(\varupalpha_{i})\]
and the \textit{entropy} of $\varupalpha$ as the expected value 
    \[H(\varupalpha) = -\sum_{i \in E}\varupmu(\varupalpha_{i})\log\varupmu(\varupalpha_{i}).\]
Consider a countable index set $S$. We  call a family $ (\varupalpha_{s})_{s\in S}$ of partitions $\varupalpha_{s} = \{\varupalpha_{s}^{i}\colon i \in E\}$ of $\Omega$ a \textit{stochastic process} with index set $S$ and state space $E$. Given a finite set $F \subseteq S$ we shall denote by $\varupalpha_{F}$ the joint partition 
    \[\bigvee_{s \in F} \varupalpha_{s} := \left\{\bigcap_{s \in F}\varupalpha_{s}^{i_{s}}\colon i_{s} \in E\right\}. \] 
    We may think of a partition $\varupalpha = \{\varupalpha^{i}\colon i \in E\}$ as well as a random variable taking values in $E$ by sending a point $\omega \in \Omega$ to the unique index $i(\omega) \in E$ with $\omega \in \varupalpha^{i(\omega)}$. This way, we may think of a family $ (\varupalpha_{s})_{s\in S}$ as above as a stochastic process in the usual sense, which justifies the naming. The partition $\varupalpha_{F}$ corresponding to a finite set $F \subseteq S$ then corresponds to a random variable taking values in $E^{F}$.
    
We shall say that a process $(\varupalpha_{s})_{s \in S}$ with finite state space $E$ satisfies the \textit{weak (mean, pointwise) equipartition property} along a sequence $(F_{n})_{n=1}^{\infty}$ of finite subsets $F_{n} \subseteq S$ if there exists a constant $h \in [0,\log |E|]$ such that 
	\[\lim_{n \to \infty}\frac{1}{|F_{n}|}~\!\!I\big(\varupalpha_{F_{n}}\big) = h\]
in probability (in $L^1$, almost surely). Obviously, both the mean and pointwise equipartition property imply the weak equipartition property. If the mean equipartition property is satisfied the constant $h$ admits a natural interpretation as an entropy rate $h$, which is given by the limit 
    \[h = \lim_{n \to \infty}\frac{H\big(\varupalpha_{F_{n}}\big)}{|F_{n}|}.\]

For stationary processes over the integers equipartition properties were first studied by Shannon in his seminal paper \cite{Shannon48}, where he proved the mean equipartition property for ergodic Markov chains. His result was extended to general ergodic procesess by McMillan in \cite{McMillan53}. Building on the work of Shannon and McMillan, Breiman established pointwise equipartition for such processes in \cite{Breiman57}. The latter result is nowadays known as the Shannon-McMillan-Breiman theorem and plays a fundamental role in entropy and information theory.

A natural question is whether the above properties can be extended to processes over more general graphs. As a first step towards generalizations, Thouvenout considered   the $d$-dimensional integer lattice $\mathbb{Z}^{d}$ in \cite{Thouvenot72}. There, the mean equipartition property was established along cubes of growing size for translation invariant and ergodic processes. This result was extended by Föllmer in \cite{Follmer73} to pointwise equipartition in the case of suitable Gibbs fields. Both results were subsequently superseded by work of Kieffer and Ornstein-Weiss in \cite{Kieffer75} and \cite{OW83}, who developed a far-reaching generalization of entropy theory for stationary processes over amenable groups. As a consequence, one obtains that for every translation invariant and ergodic process on $\mathbb{Z}^{d}$, the above mentioned results not only hold for cubes, but also for a large variety of F\o lner sequences. Building on the concept of empirically and locally converging sequences of measures, Austin studied in \cite{Aus19} asymptotic equipartition phenomena for processes over sofic groups. These results can be interpreted as relatives of Shannon-McMillan type theorems, but the underlying notion of equipartition is of a different nature.

For processes over regular trees, which are probably the most prominent model of ``multi-dimensional time'' apart from $\mathbb{Z}^{d}$, much less is known on stochastic almost everywhere entropy equipartition. A pertinent result of the above kind is due to Berger and Ye, who established the existence of entropy rates and the weak equipartition property for processes whose law is invariant and ergodic under the group of parity preserving tree automorphisms, see \cite{BY96}. Starting with \cite{Yang03}, a series of papers on strengthenings of this result for different versions of Markov chains on trees has appeared, see e.g.\@ \cite{SXZ24} and the references given therein.

\bigskip

\noindent \textbf{Main results.} In this paper we study equipartition properties for processes on regular trees from a boundary point of view. For amenable groups, this orbital approach has been pursued by  
Danilenko and Park in \cite{Dan01,DP02}. Adopting this philosophy, Nevo and Pogorzelski obtained in \cite{NP21} and \cite{NP24} pointwise and mean equipartition results for the free group along geodesics and horospheres corresponding to randomly chosen boundary points. These results are valid for processes that are invariant and ergodic under suitable groups of shift automorphisms. In this paper, we consider processes that are invariant under the full group of parity preserving tree automorphisms. Under this stronger assumption, we are able to prove  almost everywhere entropy equipartition along horospherical spheres for {\em every} choice of horosphere, see Theorem~\ref{main1} below.  
Our second main result, Theorem~\ref{main2}, states a Shannon-McMillan-Breiman theorem along metric spheres of even radius for processes which satisfy a uniform mixing condition. To best of the authors' knowledge, these are the first such results  along a deterministic sequence of spherical subsets of a non-amenable Cayley graph. 

Let us now explain the setting in more detail. Fix $d > 2$ and let $\T = (S,B)$ be the $d$-regular tree. Denote by $\partial\T$ the boundary of $\T$. Fix a root $s_{0} \in S$. Given a point $\xi \in \partial\T$, we denote by $H(\xi)$ the horosphere at $\xi$ through $s_{0}$. Writing $B_{n}$ and $S_{n}$ for the balls and spheres in $\T$ of radius $n$ with center $s_{0}$ in the path metric, we define the sets $B_{n}^{\xi} := B_{2n} \cap H(\xi)$ and $S_{n}^{\xi} := S_{2n} \cap H(\xi)$ for $n \in \mathbb{N}$, to which we refer as \textit{horospherical balls} and \textit{spheres}, respectively. Denoting by $\mathrm{Aut}_{0}(\T)$ the group of parity preserving tree automorphisms, our first main result reads as follows.

\begin{theorem} \thm \label{main1}
    Let $(\varupalpha_{s})_{s \in S}$ be a process with finite state space $E$ whose law is invariant and ergodic under the group $\mathrm{Aut}_{0}(\T)$. Then there is a constant $h \in [0,\log|E|]$ such that for every $\xi \in \partial \T$ we obtain
        \[\lim_{n \to \infty} \frac{I\big(\varupalpha_{\smash{B_{2n}^{\xi}}}\big)}{\big|B_{2n}^{\xi}\big|} = h\]
    in $L^{1}$ and almost surely.
\end{theorem}

The proof of Theorem \ref{main1} relies crucially on a result of Lubotzky-Mozes and Pemantle, cf. \cite{LM92} and \cite{Pemantle92}, stating that a process of the above type is automatically mixing in the ergodic theoretic sense. It turns out that if we impose a stronger notion of mixing we can in fact obtain equipartition along spheres. The mixing condition, to which we will refer as \textit{exponential $\varuppsi$-mixing}, requires that the stochastic dependence of events taking place in different regions $U,V \subseteq S$, as measured by the so called $\varuppsi$-coefficient, grows at most linearly in the sizes of $U$ and $V$ and decays at least exponentially in the distances of $U$ and $V$. As we shall see below, examples of processes satisfying this property arise naturally from the theory of Gibbs measures. By results of the second author in \cite{Zimmermann25+}, these examples include, for instance, the unique Gibbs states corresponding to the Ising model at high temperatures or the antiferromagnetic Potts model with sufficiently many spin states.

\begin{theorem} \thm \label{main2}
    Let $(\varupalpha_{s})_{s \in S}$ be a process with finite state space $E$ whose law is invariant under the group $\mathrm{Aut}_{0}(\T)$. If $(\varupalpha_{s})_{s \in S}$ is $\varuppsi$-mixing with exponential decay rate $\lambda > 2\log(d{-}1)$, then there exists a constant $h \in [0,\log|E|]$ such that
        \[\lim_{n \to \infty}\frac{I\big(\varupalpha_{S_{2n}}\big)}{|S_{2n}|} = h\]
    in $L^{1}$ and almost surely. 
\end{theorem} 

\medskip

\noindent \textbf{Method of proof.} As indicated already, equipartition results for processes over integer lattices arise as a special case of a far-reaching generalization of classical entropy theory to processes over countable amenable groups due to Kieffer and Ornstein-Weiss. Indeed, one can show that any stationary and ergodic process $(\varupalpha_{g})_{g \in \Gamma}$ over a countable amenable group $\Gamma$ satisfies the mean equipartition property along any F\o lner sequence $(F_{n})_{n=1}^{\infty}$ in $\Gamma$, see \cite{Kieffer75}. Moreover, if the F\o lner sequence $(F_{n})_{n=1}^{\infty}$ is nested and doubling, one obtains pointwise equipartition as well, cf. \cite{OW83}. The latter result applies for instance to all finitely generated groups of polynomial volume growth, where a suitable F\o lner sequence is given by balls in the word metric. On the other hand, it is well known that not every amenable group contains nested and doubling F\o lner sequences. Eventually, Lindenstrauss extended  in \cite{Lindenstrauss01} the Ornstein-Weiss result to so called tempered F\o lner sequences satisfying a mild growth condition. Such sequences   exist in every amenable group.

These results can be used to obtain equipartition results for processes of graphs admitting a simply transitive action of an amenable group such as the $d$-dimensional integer lattice. The $d$-regular tree $\T$ on the other hand admits a simply transitive action of the $d$-fold free product $\mathbb{F}$ of the cyclic group $\mathbb{Z}_{2}$ with itself (or of the free group of rank $d/2$ in case of an even degree). However, since these groups are not amenable, the above results cannot be directly applied to processes over the regular tree. 

However, it turns out that we can leverage the theory for amenable groups even in the setting of trees by utilizing the fact that $\mathbb{F}$ acts on $\partial\T$ in an amenable way. In the context of ergodic theorems for free groups, this approach was successfully developed by work of Bowen and Nevo, see \cite{BN13}. In the present paper, we will construct an explicit amenable group $\mathbb{G}$ acting on $\partial\T$ such that the orbits of the $\mathbb{G}$-action are contained in the orbits of the $\mathbb{F}$-action. This allows us to define a map $s\colon \mathbb{G} \times \partial\T \to S$, such that the image of $\mathbb{G}$ under the map $s(\cdot,\xi)$ equals the horosphere $H(\xi)$ for every $\xi \in \partial\T$. Now for a given process $(\varupalpha_{s})_{s \in S}$ with finite state space $E$ and a boundary point $\xi \in \partial\T$ we may consider the process $(\varupbeta^{\xi}_{g})^{~}_{g \in \mathbb{G}}$ given by
    \[\varupbeta_{g}^{\xi} := \varupalpha_{s(g,\xi)}\]
for $g \in \mathbb{G}$. If $(\varupalpha_{s})_{s \in S}$ is invariant under $\mathrm{Aut}_{0}(\T)$, it is not difficult to verify that the processes $(\varupbeta^{\xi}_{g})^{~}_{g \in \mathbb{G}}$ are stationary. Moreover, we will show that the processes $(\varupbeta^{\xi}_{g})^{~}_{g \in \mathbb{G}}$ are ergodic, in fact mixing, whenever the process $(\varupalpha_{s})_{s \in S}$ is ergodic under $\mathrm{Aut}_{0}(\T)$. To this end we will make use of the fact that every ergodic action of the group $\mathrm{Aut}_{0}(\T)$ on a standard probability space is automatically mixing. This follows from a theorem of Lubotzky and Mozes stating that $\mathrm{Aut}_{0}(\T)$ has the Howe-Moore property, cf. \cite{LM92}. In the specific case of processes over regular trees a similar result was shown by Pemantle in \cite{Pemantle92}. 

We are therefore in a situation, where we may apply the theorems of Ornstein-Weiss and Lindenstrauss to obtain that for any $\xi \in \partial \T$ the process $(\varupbeta_{g}^{\xi})^{~}_{g \in \mathbb{G}}$ satisfies the mean and pointwise equipartition property along any suitable F\o lner sequence $(\smash{F^{\xi}_{n}})_{n=1}^{\infty}$. However, by definition of the processes $(\varupbeta_{g}^{\xi})^{~}_{g \in \mathbb{G}}$ this is equivalent to the statement that the process $(\varupalpha_{s})_{s \in S}$ satisfies the mean and pointwise equipartition property along the sequence$(\smash{\widehat{F}_{n}^{\xi}})_{n=1}^{\infty}$ of sets given by
    \[\widehat{F}_{n}^{\xi} := \{s(g,\xi)^{-1}\colon g \in F_{n}\}\]
for $n \in \mathbb{N}$. Choosing appropriate F\o lner sequences will then lead to the proof of Theorem~\ref{main1}. To prove Theorem~\ref{main2}, we consider the canonical Patterson-Sullivan measure $\varupnu$ on the boundary $\partial\T$. Using a suitable maximal inequality we may integrate over this probability measure to obtain a constant $h$ such that
    \[\int_{\partial\T} \frac{I\big(\varupalpha^{~}_{\smash{\widehat{F}_{n}^{\xi}}}\big)}{\big|\widehat{F}_{n}^{\xi}\big|}~\mathrm{d}\varupnu(\xi) \to h\]
$\varupmu$-almost surely. Now we will choose the F\o lner sequences $(F_{n}^{\xi})_{n=1}^{\infty}$ in such a way that for every $n \in \mathbb{N}$ we obtain sets $\widehat{F}_{n,1}, \ldots, \widehat{F}_{n,r_{n}}$ with
    \[S_{2n} = F_{n,1} \dot{\cup} \ldots \dot{\cup} F_{n,r_{n}}\]
and    
    \[\int_{\partial\T} \frac{I\big(\varupalpha^{~}_{\smash{\widehat{F}_{n}^{\xi}}}\big)}{\big|F_{n}^{\xi}\big|}~\mathrm{d}\varupnu(\xi) = \frac{1}{|S_{2n}|}\sum_{j=1}^{r_{n}}I\big(\varupalpha_{\smash{\widehat{F}_{n,j}}}\big).\]
Moreover, we will make sure that the distance of the sets $\widehat{F}_{n,i}$ and $\widehat{F}_{n,j}$ grows linearly in $n$ for $i \neq j$. These properties, together with the mixing condition, ensure that for all $i \neq j$ and every $\omega \in \Omega$ the atoms $\smash{{\varupalpha_{\vphantom{\widehat{F}}}^{~\!\!i(\omega)}}}_{\!\!\!\!\!\!\!\!\!\!\!\!\widehat{F}_{n,i}}$ and $\smash{{\varupalpha_{\vphantom{\widehat{F}}}^{~\!\!i(\omega)}}}_{\!\!\!\!\!\!\!\!\!\!\!\!\widehat{F}_{n,j}}$ become more and more stochastically independent as $n$ goes to $\infty$. Thus we obtain 
    \[\left|\frac{1}{|S_{2n}|}\sum_{j=1}^{r_{n}}I\big(\varupalpha_{\widehat{F}_{n,j}}\big) - \frac{1}{|S_{2n}|}I\big(\varupalpha_{S_{2n}}\big)\right| \to 0\]
$\varupmu$-almost surely for $n \to \infty$. This yields the claimed  almost sure equipartition property. From this and the above mentioned maximal inequality mean equipartition is then  obtained via standard methods.

\medskip

\noindent \textbf{Organization of the paper}. In Chapter~2, we provide the relevant preliminaries on graph automorphisms, group actions, random processes and Borel equivalence relations. We study the synchronous tail relation over the boundary of regular trees in Chapter~3, representing  the equivalence classes as orbits of a locally finite abelian group. Chapter~4 is devoted to the proof of almost everywhere equipartition along horospherical spheres. Under the assumption of exponential $\psi$-mixing we prove almost everywhere equipartition along metric spheres in Chapter~5.

\medskip

\noindent \textbf{Acknowledgements}. The authors wish to express their gratitude to Amos Nevo and Anush Tserunyan 
for many insightful and stimulating conversations. FP thanks Technion Haifa and IIAS Jerusalem for excellent working conditions during research stays in 2023 and 2025. The results of this work are part of the doctoral project of EZ.   
The authors gratefully acknowledge support by grant I-1485/304.6-2019
from the German Israeli Foundation for Scientific Research and Development (GIF). EZ was supported by a doctoral scholarship from the German Academic Scholarship Foundation (Studienstiftung des deutschen Volkes).

\section{Preliminaries}

\noindent \textbf{Graphs and automorphisms.} By a \textit{graph} we mean a tuple $G = (S,B)$ consisting of a countable set $S$ of {\em sites} or {\em vertices} and a symmetric set $B \subseteq S \times S$ of {\em bonds} or {\em edges}. For any $s \in S$ we define the \textit{degree} $\deg(s)$ of $s$ as the number of $t \in S$ such that $(s,t) \in B$. We call ${G}$ \textit{regular} of degree $d$ if $\deg(s) = d$ for all $s \in S$. A \textit{path} in ${G}$ is given by a sequence $s_{1},\ldots,s_{m}$ such that $(s_{i},s_{i+1}) \in B$ for $i = 1,\ldots,m{-}1$. We call $m$ the \textit{length} of the path. If $s_{1} = s_{m}$ the path is called a \textit{cycle}. The graph $G$ is called \textit{connected} if for all $s,t \in S$ there is a path $s_{1},\ldots,s_{m}$ with $s=s_{1}$ and $t = s_{m}$. In this case we define a metric $d$ on $S$ by assigning to every pair $s,t \in S$ the length of a minimal path connecting $s$ and $t$. Given a site $s \in S$ we shall denote by 
	\[B_{n}(s) := \{t \in S\colon d(s,t) \leq n\}\]
the \textit{ball} of radius $n$ and by 
	\[S_{n}(s) := \{t \in S\colon d(s,t) = n\}\]
the \textit{sphere} of radius $n$ respectively. For every graph appearing in this work, we consider its vertex set as a metric space with the path metric. 
We shall call a graph \textit{acyclic} if it contains no cycles. A connected and acyclic graph is called a \textit{tree}.

Let us assume that $G$ is regular and connected. By a \textit{graph automorphism} we mean a bijection $\varupphi\colon S \to S$ satisfying $(\varupphi(s),\varupphi(t)) \in B$ if and only if $(s,t) \in B$. The collection of graph automorphisms constitutes a group, which we shall denote by $\mathrm{Aut}(G)$. Equipping $\mathrm{Aut}(G)$ with the subspace topology inherited from the product topology on the space of all maps from $S$ to itself, we obtain a locally compact and second countable Polish group. Fixing a vertex $s_{0} \in S$ we may define a complete metric $d$ on $\mathrm{Aut}(G)$ generating the topology by setting
	\[d(\varupphi,\varuppsi) := 2^{-n(\upphi,\uppsi)}\]
for $\varupphi,\varuppsi \in \mathrm{Aut}(G)$, where $n(\varupphi,\varuppsi)$ denotes the maximal radius $n \in \mathbb{N}_{0}$ such that $\varupphi(s) = \varuppsi(s)$ for all $s \in B_{n}(s_{0})$. 

Let $\Gamma = \langle \Sigma \rangle$ be a finitely generated group with a symmetric set of generators $\Sigma$. By the \textit{(right) Cayley graph} of $\Gamma$ with respect to $\Sigma$ we mean the graph $G$ with vertex set $S = \Gamma$ and edge set $B$ consisting all pairs $(g,h)$ such that $h = gs$ for some $s \in \Sigma$. Cayley graphs are regular and connected and the path metric corresponds to the word metric of $\Gamma$ with respect to $\Sigma$. In this case every $g \in \Gamma$ induces an automorphism $\varuptheta_{g} \in \mathrm{Aut}(G)$ given by $\varuptheta_{g}(h) := gh$ for $h \in \Gamma$. This way we may consider $\Gamma$ as a (discrete) subgroup of $\mathrm{Aut}(G)$.

\medskip

\noindent \textbf{Group actions.} Let $(X,\mathscr{B})$ be a standard Borel space and $G$ be a locally compact second countable Hausdorff group. Let $\mathrm{Aut}(X)$ denote the group of Borel automorphisms of $X$. By a \textit{measurable action} we mean a group homomorphism $\varuppi\colon G \to \mathrm{Aut}(X)$ such that the mapping $(g,x) \mapsto \varuppi(g)x$ is measurable. If $\varuppi$ is clear from the context, we will denote the automorphism $\varuppi(g)$ just by $g$. We shall call a Borel probability measure $\varupmu$ on $X$ \textit{invariant} under the action of $G$ on $X$ if $\varupmu(gB) = \varupmu(B)$ for all $B \in \mathscr{B}$. A set $B \in \mathscr{B}$ is called \textit{invariant} if $gB = B$ for every $g \in G$. The action of $G$ on $X$ is said to be \textit{ergodic} if $\varupmu(B) \in \{0,1\}$ for any invariant set $B \in \mathscr{B}$. Given a sequence $(g_{n})_{n=1}^{\infty}$ in $G$ we shall write $g_{n} \to \infty$ to indicate that for every compact set $K \subseteq G$ there is some index $n_{0} \in \mathbb{N}$ such that for all $n \geq n_{0}$ one has $g_{n} \in K^{c}$. The action of $G$ on $X$ is said to be \textit{mixing} if for every sequence $(g_{n})_{n=1}^{\infty}$ in $G$ with $g_{n} \to \infty$ and for all sets $B,C \in \mathscr{B}$ we obtain 
    \[\lim_{n \to \infty} \varupmu(A \cap g_{n}B) = \varupmu(A)\varupmu(B).\]

Let $\Lambda$ be a countable group with neutral element $e$ acting on a standard Borel space $(X,\mathscr{B})$ in a measurable fashion. Let $\Gamma$ be a further countable group with neutral element $e$. By a \textit{Borel cocycle} we mean a measurable map $\varupvarphi\colon \Lambda \times X \to \Gamma$ such that $\varupvarphi(1,x) = e$ and 
    \[\varupvarphi(gh,x) = \varupvarphi(g,hx)\varupvarphi(h,x)\]
for every $x \in X$ and all $g,h \in \Lambda$. 

\medskip

\noindent \textbf{Stochastic processes.} Let $(\Omega,\mathscr{A},\varupmu)$ be a standard probability space. By a \textit{(finite) measurable partition} of $\Omega$ we mean a finite collection $\varupalpha = \{\varupalpha^{i}\colon i \in E\}$ consisting of disjoint measurable sets $\varupalpha^{i}$ whose union equals $\Omega$. Given two such partitions $\varupalpha = \{\varupalpha^{i}\colon i \in E\}$ and $\varupbeta = \{\varupbeta^{~\!j}\colon j \in F\}$ we define their \textit{join} by 
    \[\varupalpha \vee \varupbeta := \left\{\varupalpha^{i} \cap \varupbeta^{~\!j}\colon (i,j) \in E \times F\right\}.\]
Obviously $\varupalpha \vee \varupbeta$ is again a finite measurable partition. Consider a countable set $S$ and fix a finite set $E$. By a \textit{stochastic process} on $S$ taking values in $E$ we mean a family $(\varupalpha_{s})_{s \in S}$ consisting of finite measurable partitions $\varupalpha_{s} = \{\varupalpha_{s}^{i}\colon i \in E\}$ of $\Omega$. Given a finite set $F \subseteq S$ we shall set
    \[\varupalpha^{~}_{F} := \bigvee_{s \in F} \varupalpha_{s}.\]
Equipping the set $E^{S}$ with the product topology we obtain a compact Polish space whose Borel $\sigma$-algebra we shall denote by $\mathscr{B}$. As indicated earlier, we may think of the partitions $\varupalpha_{s}$ also as random variables $\varupalpha_{s}\colon \Omega \to E$ sending a point $\omega \in \Omega$ to the unique $i(\omega) \in E$ such that $\omega \in \varupalpha_{s}^{i(\omega)}$. This way we may identify the partitions $\varupalpha_{F}$ with random variables taking values in $E^{F}$ and the full process with a random variable taking values in $E^{S}$. The \textit{law} of of the process $(\varupalpha_{s})_{s \in S}$ is given by the push-forward measure $\mathrm{P}$ of $\varupmu$ under the the latter map. The arising standard probability space $(E^{S},\mathscr{B},\mathrm{P})$ is called the \textit{configuration space} of the process $(\varupalpha_{s})_{s \in S}$. We shall say that the process $(\varupalpha_{s})_{s \in S}$ is given in \textit{canonical form} if $(\Omega,\mathscr{A},\varupmu)$ equals the configuration space and $\varupalpha_{s}$ equals the the partition $\{\varuppi^{-1}_{s}(i)\colon i \in E\}$, where $\varuppi_{s}\colon \Omega \to E$ denotes the projection to the coordinate $s$ for $s \in S$. 

Let $G$ be a locally compact second countable Hausdorff group acting on $E^{S}$ in a measurable fashion. We will be interested in situations where the law $\mathrm{P}$ of $(\varupalpha_{s})_{s \in S}$ is invariant (ergodic, mixing) with respect to the action of $G$. A specific instance of the foregoing setting is obtained when $G$ admits a measurable action on $S$ via permutations, in which case we may consider the induced (measurable) action of $G$ on $E^{S}$ given by
    \[g(x_{s})_{s \in S} := (x_{g^{-1}s})_{s \in S}\]
for $g \in G$ and $x \in X$. In this case invariance under the action of $G$ amounts to the assumption that for every every $g \in G$ the laws of the processes $(\varupalpha_{s})_{s \in S}$ and $(\varupalpha_{g^{-1}s})_{s \in S}$ are the same. 
Moreover, mixing corresponds to the condition that for every sequence $(g_{n})_{n=1}^{\infty}$ with $g_{n} \to \infty$, all $U,V \subseteq S$ finite and every $\varepsilon > 0$ there is some $n_{0} \in \mathbb{N}$ such that we have
    \[\big|\varupmu\big(\varupalpha_{g_{n}V}^{~\!\!i} \cap \varupalpha_{U}^{~\!\!j}\big) - \varupmu\big(\varupalpha_{V}^{~\!\!i}\big)\varupmu\big(\varupalpha_{U}^{~\!\!j}\big)\big| < \varepsilon\]
for all $i,j$ and $n \geq n_{0}$. A special case of the above setting arises when $S$ is the set of sites of a regular graph $G$ and we consider the action of the automorphism group $\mathrm{Aut}(G)$ on $S$. Another special case is obtained when $S$ is a countable group $\Gamma$ and we consider the action of $\Gamma$ on itself by left multiplication. In the latter case we shall call a process $(\varupalpha_{g})_{g \in \Gamma}$ whose law is invariant (ergodic, mixing) under the action of $\Gamma$ \textit{stationary (ergodic, mixing)}. 

\medskip

\noindent \textbf{Measurable equivalence relations.} Let $(X,\mathscr{B})$ be a standard Borel space. By a \textit{countable Borel equivalence relation} we mean a Borel set $E \subseteq X \times X$, which is an equivalence relation whose equivalence classes are countable. We call a Borel automorphism $\varupphi\colon X \to X$ an \textit{inner automorphism} of $E$ if $\textrm{graph}(\varupphi) \subseteq E$. Let $\Gamma$ be a countable group acting on $X$ in a measurable fashion. Then we call an equivalence relation $E$ the \textit{orbit relation} of the $\Gamma$-action if it consists of pairs of pairs $(x,gx)$ for $x \in X$ and $g \in \Gamma$. In this case we will also say that $E$ is \textit{generated} by the action of $\Gamma$. Orbit relations are always countable Borel.

\section{Regular trees and boundary actions} \label{Section_Boundary}

Fix $d > 2$. Let $\mathbb{F}$ be the $d$-fold free product $\mathbb{Z}_{2} * \ldots * \mathbb{Z}_{2}$ of the cyclic group $\mathbb{Z}_{2}$ with itself. The $d$-regular tree $\T = (S,B)$ arises as the Cayley graph of $\mathbb{F}$ with respect to the standard set of generators $\Sigma := \{a_{1},\ldots,a_{d}\}$. Thus we may identify the elements of $S$ with finite words $s_{1}\ldots s_{n}$ over the alphabet $\Sigma$ satisfying $s_{i} \neq s_{i+1}$ for $i = 1,\ldots,n$. The boundary $\partial\T$ of $\T$ may then be defined as the set of all one-sided infinite sequences $\xi$ over $\Sigma$ satisfying $\xi_{n} \neq \xi_{n+1}$ for $n \in \mathbb{N}$. (If $d$ is even, we could as well choose $\mathbb{F}$ to be the free group of rank $d$ with standard set of generators $\Sigma = \{a_1^{\pm 1}, \dots, \smash{a_{d/2}^{\pm 1}}\}$, in which case the forbidden relations are given by $s_{i+1} = s_i^{-1}$ and $\xi_{n+1} = \xi_n^{-1}$.) We may equip $\partial\T$ with the subspace topology inherited from the product topology on $\Sigma^{\mathbb{N}}$. As a closed subset of $\Sigma^{\mathbb{N}}$ it defines a compact Polish space. It is homeo\-morphic to the Gromov boundary associated to $\mathbb{F}$ as a hyperbolic group as well as to the space of ends associated to the tree $\T$ as a graph. The \textit{Patterson-Sullivan measure} $\varupnu$ on $\partial\T$ is given by the unique Markov measure corresponding to the initial probability vector $p = (p_{i})_{i=1}^{d}$ and the transition matrix $P = (p_{ij})_{ij = 1}^{d}$ defined by $p_{i} = d^{-1}$ for $i = 1,\ldots,d$ and 
    \[p_{ij} := \frac{1}{d{-}1}(1{-}\delta_{ij})\]
for $i,j = 1,\ldots,d$. (If $\mathbb{F}$ is a free group of even rank $d$ with standard set of generators, then $p_{ij} = 1/(d-1)$ as long as $s_{i+1} \neq s_i^{-1}$.)

The group $\mathbb{F}$ acts on $\partial\T$ via homeomorphisms by concatenation and cancellation. We will denote the respective orbit relation by $R$. Let $R_{0}$ be the equivalence relation on $\partial\T$ consisting of all pairs $(\xi,\zeta)$ that differ only in finitely many coordinates. It is not difficult to see that $R_{0}$ is a Borel subrelation of $R$. Denoting for $n \geq 0$ by $R_{0}^{n}$ the equivalence relation consisting of all pairs $(\xi,\zeta)$ that differ only in the first $n$ entries (where we define $R_{0}^{0}$ to be the identity relation) we obtain
    \[R_{0} = \bigcup_{n=0}^{\infty} R_{0}^{n}.\]
Noting that the equivalence relations $R_{0}^{n}$ are Borel and have finite classes this shows that $R_{0}$ is hyperfinite. In particular, by a classical result of Slaman and Steel, see \cite{SlSt88}, $R_{0}$ is the orbit relation of a Borel action of $\mathbb{Z}$ on $\partial\T$. However, defining such a $\mathbb{Z}$-action in an explicit fashion is notoriously difficult. On the other hand, considering 
an alternative amenable group, we may explicitly construct Borel actions generating $R_{0}$ as follows.

Let us define for every $n \in \mathbb{N}_{0}$ and every $R_{0}^{n}$-class $C$ a linear ordering $\prec_{C}$ on $C$ by the following inductive procedure. In the induction basis there is only one choice of an ordering since the equivalence classes of $R_{0}^{0}$ are singletons. For the induction step fix $n \geq 0$ and let $C$ be some equivalence class of $R_{0}^{n+1}$. Noting that every $R_{0}^{n}$-class consists of sequences with a common $n$-th entry we may order the $R_{0}^{n}$-classes $C_{1},\ldots,C_{d-1}$ contained in  $C$ according to the $n$-th entries of their elements. Fix distinct $\xi,\zeta \in C$ with $\xi \in C_{i}$ and $\zeta \in C_{j}$ and assume that $\xi$ is the $k$-th element in $C_{i}$ and $\zeta$ is the $l$-th element in $C_{j}$. Then we stipulate that $\xi \prec_{C} \zeta$ if and only if either $k < l$ or ($k = l$ and $i < j$). It is straightforward to check that $\prec_{C}$ defines indeed a linear ordering. 

Now for every $n \in \mathbb{N}_{0}$ we shall consider the map $\varuppsi_{n}\colon \partial \T \to \partial \T$ acting on an equivalence class $C$ of $R_{0}^{n}$ as a cyclic permutation by sending the maximal element to the minimal element and every other element to its successor with respect to $\prec_{C}$. Let $\Psi_{n}$ denote the group of maps generated by $\varuppsi_{n}$ and let $\Psi$ be the smallest group containing all of these groups.

\begin{theorem} \label{Proposition} \proposition
The maps $\varuppsi_{n}$ are homeomorphisms of $\partial\T$ and inner automorphisms of $R_{0}^{n}$. Furthermore, the groups $\Psi_{n}$ are cyclic groups of order $(d{-}1)^{n}$ satisfying $R_{0}^{n} = \textnormal{graph}(\Psi_{n})$ and $\Psi_{n} \subseteq \Psi_{n+1}$ for every $n \in \mathbb{N}_{0}$. In particular, we have $\Psi = \bigcup_{n=0}^{\infty} \Psi_{n}$ and $R_{0} = \textnormal{graph}(\Psi)$.
\end{theorem}

\textit{Proof:} Fix $n \in \mathbb{N}_{0}$. Obviously the map $\varuppsi_{n}$ is bijective and satisfies $\text{graph}(\varuppsi_{n}) \subseteq R_{0}^{n}$. The latter implies that $\varuppsi_{n}$ fixes all except maybe the first $n$ letters of a sequence $\xi$. It is furthermore easy to check that $\varuppsi_{n}(\xi)$ depends only on the first $n$ letters of $\xi$. Both together shows that $\varuppsi_{n}$ as well as $\varuppsi_{n}^{-1}$ map cylinder sets to cylinder sets, so both are continuous and in particular Borel. In sum this implies that $\varuppsi_{n}$ is a homeomorphism of $\partial\T$ as well as an inner automorphism of $R_{0}^{n}$. The fact that $\varuppsi_{n}$ is a cyclic group of order $(d{-}1)^{n}$ follows from the fact that the $R_{0}^{n}$-classes have precisely $(d{-}1)^{n}$ elements. Furthermore, the construction of $\varuppsi_{n}$ implies that $\text{graph}(\Psi_{n}) = R_{0}^{n}$.

To verify that $\Psi_{n} \subseteq \Psi_{n+1}$ fix $\xi \in \partial\T$, let $C$ be the $R_{0}^{n+1}$-class of $\xi$ and $C_{1},\ldots,C_{d-1}$ be the ordered $R_{0}^{n}$-classes contained in $C$. Let $C_{i}$ the the $R_{0}^{n}$-class containing $\xi$ and set $\zeta := \varuppsi_{n}(\xi)$. Then we may distinguish two cases. We  consider first the case that $\xi$ is not the maximal element in $C_{i}$, so $\zeta$ is the successor of $\xi$ in $C_{i}$. Let us assume that $\xi$ is the $j$-th element and $\zeta$ the $(j{+}1)$-th element in $C_{i}$. Then there are precisely $d{-}2$ elements between $\xi$ and $\zeta$ in $C$, namely the $j$-th elements of the classes $C_{i+1}, \ldots , C_{d-1}$ and the $(j{+}1)$-th elements of the classes $C_{1},\ldots, C_{i-1}$. Thus we have $\varuppsi_{n+1}^{d-1}(\xi) = \zeta = \varuppsi_{n}(\xi)$. Now consider the case that $\xi$ is the maximal element in $C_{i}$, so $\zeta$ is the minimal element in $C_{i}$. Then there are precisely $d{-}i{-}1$ elements of $C$ above $\xi$, namely the maximal elements of the classes $C_{i+1}, \ldots , C_{d-1}$, and $i{-}1$ elements of $C$ below $\zeta$, namely the minimal elements of the classes $C_{1},\ldots, C_{i-1}$. Thus we have again $\varuppsi_{n+1}^{d-1}(\xi) = \zeta = \varuppsi_{n}(\xi)$. Both together shows $\varuppsi_{n} = \varuppsi_{n+1}^{d-1}$, which implies $\Psi_{n} \subseteq \Psi_{n+1}$. This gives $\Psi = \bigcup_{n=0}^{\infty}\Psi_{n}$ and thus $R_{0} = \text{graph}(\Psi)$. \hfill $\Box$ \\

As an abstract group $\Psi$ has the structure of a direct limit of cyclic groups. To see this we consider for every $n \in \mathbb{N}_{0}$ the root of unity
	\[g_{n} := \exp\bigg(\frac{2\pi i}{(d{-}1)^{n}}\bigg)\]
 and identify the cyclic group of order $(d{-}1)^{n}$ with the multiplicative subgroup of $\mathbb{S}^{1}$ given by
	\[\mathbb{G}_{n} := \langle g_{n} \rangle = \Big\{g_{n}^{m}\colon 0 \leq m \leq (d{-}1)^{n} - 1\Big\}.\]
Then we obtain $\mathbb{G}_{0} \subseteq \mathbb{G}_{1} \subseteq \mathbb{G}_{2} \subseteq \ldots$, so we may identify the direct limit of these cyclic groups with the union
	\[\mathbb{G} = \bigcup_{n=0}^{\infty}\mathbb{G}_{n}.\]
Note that as a union of finite (and thus amenable) groups the group $\mathbb{G}$ is itself amenable. Moreover, as a subgroup of the abelian group $\mathbb{S}^{1}$ the group $\mathbb{G}$ is itself abelian. We will show that $\varuppsi$ is isomorphic to $\mathbb{G}$ as an abstract group. To this end consider the map $\Phi\colon \Psi \to \mathbb{G}$ defined by 
     \[\Phi(\varuppsi_{n}^{m}) = g_{n}^{m}\]
for $n \in \mathbb{N}_{0}$ and $0 \leq m \leq (d{-}1)^{n} -1$. 
\begin{theorem}
    \proposition The map $\Phi$ defines an isomorphism between the groups $\Psi$ and $\mathbb{G}$ satisfying $\Phi(\Psi_{n}) = \mathbb{G}_{n}$ for every $n \in \mathbb{N}_{0}$.
\end{theorem}

\textit{Proof:} As in the proof of Proposition \ref{Proposition} we obtain that $\varuppsi_{n} = \varuppsi_{n+1}^{d-1}$ for every $n \in \mathbb{N}_{0}$, which inductively implies  
    \[\varuppsi_{n} = \varuppsi_{n+j}^{(d-1)^{j}}\]
for all $j \in \mathbb{N}$. Thus, fixing $n,\ell \in \mathbb{N}_{0}$ and $m,k \in \mathbb{Z}$ we have
    \begin{align} \label{eq0}
        \varuppsi_{n}^{m} \circ \varuppsi_{\ell}^{k} = \varuppsi_{n+\ell}^{m(d-1)^{\ell}} \circ \varuppsi_{n+\ell}^{k(d-1)^{n}} = \varuppsi_{n+\ell}^{m(d-1)^{\ell}+k(d-1)^{n}}.
    \end{align}
Accordingly we obtain that $\varuppsi_{n}^{m} = \varuppsi_{\ell}^{k}$ if and only if $m(d{-}1)^{-n} = k(d{-}1)^{-\ell}$. The latter condition is satisfied if and only 
    \[g_{n}^{m} = \exp\bigg( \frac{2\pi i m}{(d{-}1)^{n}}\bigg) = \exp\bigg( \frac{2\pi i k}{(d{-}1)^{\ell}}\bigg) = g_{\ell}^{k}.\]
This shows that $\Phi$ is well defined and injective. It is also obvious that $\Phi$ is surjective. Using again $(\ref{eq0})$
we obtain furthermore
    \begin{align*}
    \Phi(\varuppsi_{n}^{m} \circ \varuppsi_{\ell}^{k}) &= \Phi\Big(\varuppsi_{n+\ell}^{m(d-1)^{\ell}+k(d-1)^{n}}\Big) \\ &= \exp\bigg( \frac{2\pi i (m(d{-}1)^{\ell}+k(d{-}1)^{n})}{(d{-}1)^{n+ \ell}}\bigg) \\
    &= \exp\bigg( \frac{2\pi i m}{(d{-}1)^{n}}\bigg) \cdot \exp\bigg( \frac{2\pi i k}{(d{-}1)^{\ell}}\bigg) \\
    &= \vphantom{\int} \Phi(\varuppsi_{n}^{m}) \cdot \Phi(\varuppsi_{\ell}^{k}).
    \end{align*}
This shows that $\Phi$ is a group isomorphism. The fact that $\Phi(\Psi_{n}) = \mathbb{G}_{n}$ follows from the definition. \hfill $\Box$ 

\medskip

We may summarize the above results by saying that the group $\mathbb{G}$ acts freely on $\partial \T$ via homeomorphisms and the corresponding orbit relation equals $R_{0}$. Moreover, for every $n \in \mathbb{N}_{0}$ the orbits of the induced action of the subgroup $\mathbb{G}_{n}$ coincide with the classes of $R_{0}^{n}$. 

\section{Equipartition along horospheres}
    This section is devoted to the proof of our first main result, Theorem \ref{main1}. We shall adopt the setting of the foregoing section. Given a pair $(\xi,\zeta) \in R_{0}$ we shall denote by $u(\xi,\zeta)$ the shortest element $u \in \mathbb{F}$ satisfying $\zeta = u^{-1}\xi$. Choosing $n \in \mathbb{N}_{0}$ with $(\xi,\zeta) \in R_{0}^{n}$ this element takes the form
        \[u(\xi,\zeta) = \xi_{1}\ldots\xi_{n}\zeta_{n}^{-1}\ldots\zeta_{n}^{-1}.\]
    Based on this we may define the map $\varupvarphi\colon \mathbb{G} \times \partial\T \to \mathbb{F}$ by setting
        \[\varupvarphi(g,\xi) := u(g\xi,\xi)\]
    for $g \in \mathbb{G}$ and $\xi \in \partial\T$. It is not difficult to check that $\varupphi$ is a Borel cocycle. Let us denote by $s(g,\xi)$ the site in $S$ corresponding to the element $u(\xi, g^{-1}\xi) = \varupvarphi(g^{-1},\xi)^{-1}$ in $\mathbb{F}$. It is not hard to see that for every $\xi \in \partial\T$ the map $s(\cdot,\xi)$ is injective. The set
        \[H(\xi) := \{s(g,\xi)\colon g \in \mathbb{G}\}\]
    is called the \textit{horosphere} at $\xi$ through $e$. Intuitively it consists of all sites $s \in S$ having ``equal distance'' to $\xi$ than $e$. 

    It is well known that the tree $\T$ is bipartite, i.e.\@ it admits a coloring of its vertices with two colors such that no adjacent vertices get the same color. A particular important closed subgroup of the automorphism group $\mathrm{Aut}(\T)$ is given by the group $\mathrm{Aut}_{0}(\T)$ of automorphisms preserving this coloring. We will call such automorphisms {\em parity preserving}. $\mathrm{Aut}_{0}(\T)$ is a closed subgroup of $\mathrm{Aut}(\T)$. It turns out that for any two horospheres we can find a parity preserving automorphism mapping one to the other while preserving the respective ``parametrizations'' by $\mathbb{G}$. To show this we will derive a series of lemmas. 

    We shall start with fixing some notation and terminology: we call two elements $u,v \in S$ \textit{compatible} if the last letter of $u$ does not equal the first letter of $v$. (If $\mathbb{F}$ is a free group of even rank, then the condition is that the inverse of the last letter of $u$ does not coincide with the first letter of $v$.) 
    Given an element $u \in S$ we shall denote by $S_{u}$ the {\em subtree under $u$}, which consists of all elements in $S$ of the form $uv$ for some $v \in S$ compatible with $u$. Furthermore, given $k \in \mathbb{N}_{0}$, we define the \textit{$k$-th level} $L_{u}^{k}$ of $S_{u}$ as the set of all elements of $S_{u}$ of length $|u| + k$. The first lemma shows that we can always find tree automorphisms that ``flip'' two subtrees arising from two children of the same vertex. We will refer to automorphisms of this type accordingly as \textit{flip automorphisms} in what follows.

\begin{theorem} \lemma   \label{flip}
	For every $u \in S$ and all $a,b \in \Sigma$ compatible with $u$ there is an automorphism $\varupphi \in \mathrm{Aut}_{0}(\T)$ mapping $S_{ua}$ to $S_{ub}$ and leaving every vertex outside of $S_{ua} \cup S_{ub}$ fixed.
\end{theorem}

\textit{Proof:} Fix $u \in S$ and $a,b \in \Sigma$ compatible with $u$. We define the automorphism $\varupphi$ as follows. For any $k \geq 0$ we order the $(d{-}1)^{k}$ elements in $L_{ua}^{k}$ and $L_{ub}^{k}$ in lexicographical order. Then we stipulate that $\varupphi$ sends the $i$-th element of $L_{ua}^{k}$ to the $i$-th element of $L_{ub}^{k}$ for $i = 1,\ldots (d{-}1)^{k}$ and $k \geq 0$ and leaves every vertex outside of $S_{ua} \cup S_{ub}$ fixed. The resulting map $\varupphi$ is obviously a bijection and satisfies the statement of the lemma. To show that $\varupphi$ preserves the tree structure it suffices to show that for any $v \in S_{ua}$ and any $c \in \Sigma$ compatible with $v$ there is some $d \in \Sigma$ compatible with $\varupphi(v)$ such that $\varupphi(vc) = \varupphi(v)d$. To see this let $k$ be the level of $v$ in $S_{ua}$ and assume that $v$ is the $i$-th element in $L_{ua}^{k}$. Then the elements of the form $vc$ have an index between $(i{-}1)(d{-}1){+}1$ and $i(d{-}1)$ in the lexicographical ordering of $L_{ua}^{k+1}$. Thus by construction $\varupphi(v)$ has index $i$ in $L_{ub}^{k}$ and $\varupphi(vc)$ has index between $(i{-}1)(d{-}1){+}1$ and $i(d{-}1)$ in $L_{ub}^{k+1}$. But this means that there is some $d \in \Sigma$ compatible with $\varupphi(v)$ such that $\varupphi(vc) = \varupphi(v)d$. The observation that $\varupphi \in \mathrm{Aut}_{0}(\T)$  finishes the proof. \hfill $\Box$ \\

The next lemma shows that for any two geodesics in $\T$ we can find a tree automorphism mapping one to the other, while preserving the distance to the root.

\begin{theorem} \lemma \label{automorphism_geodesic}
	For any two boundary points $\xi,\zeta \in \partial\T$ there is an automorphism $\varupphi \in \mathrm{Aut}_{0}(\T)$ such that $\varupphi(e) = e$ and $\varupphi(\xi_{1}\ldots\xi_{n}) = \zeta_{1}\ldots\zeta_{n}$ for every $n \in \mathbb{N}$.
\end{theorem}

\textit{Proof:} For $n \in \mathbb{N}$ set $u^{n} := \xi_{1}\ldots\xi_{n}$ and $v^{n} := \zeta_{1}\ldots \zeta_{n}$. We will show that there exists a Cauchy sequence $(\varupphi_{k})_{k =1}^{\infty}$ of automorphisms $\varupphi_{k} \in \mathrm{Aut}_{0}(\T)$ such that for every $k \in \mathbb{N}$ we have 
	\[\varupphi_{k}(u^{i}) = v^{i}\]
for $1 \leq i \leq k$. Using the fact that $\mathrm{Aut}(\T)$ is a complete space and $\mathrm{Aut}_{0}(\T)$ is a closed subgroup we obtain a limit $\varupphi \in \mathrm{Aut}_{0}(\T)$. To see that $\varupphi$ satisfies the claim, fix $n \in \mathbb{N}$. Then for all $k \geq n$ large enough $\varupphi_{k}$ coincides with $\varupphi$ on the ball $B_{n}$ of radius $n$. Noting that $u^{n} \in B_{n}$ this gives $\varupphi(u^{n}) = \varupphi_{k}(u^{n}) = v^{n}$, which shows the claim. 

So it remains to construct the Cauchy sequence $(\varupphi_{k})_{k=1}^{\infty}$, which will be done inductively as follows. For $k = 1$ Lemma \ref{flip} provides us with a flip automorphism $\varupphi_{1} \in \mathrm{Aut}_{0}(\T)$ satisfying $\varupphi_{1}(e) = e$ and $\varupphi_{1}(u^{1}) = v^{1}$. Now assume that we defined $\varupphi_{k} \in \mathrm{Aut}_{0}(\T)$ such that for $1 \leq i \leq k$ we have $\varupphi_{k}(u^{i}) = v^{i}$. Then we have in particular $\varupphi_{k}(u^{k}) = v^{k}$. Noting that the vertices $\varupphi_{k}(u^{k+1})$ and $\varupphi_{k}(u^{k})$ are neighbours, the latter implies that the vertices $\varupphi_{k}(u^{k+1})$ and $v^{k}$ are neighbours. Since $\varupphi_{k}(u^{k-1}) = v^{k-1}$ this is only possible if there is some $a \in \Sigma$ compatible with $v^{k}$ such that $\varupphi_{k}(u^{k+1}) = v^{k}a$. By Lemma \ref{flip} there is a flip automorphism $\varuppsi \in \mathrm{Aut}_{0}(\T)$ satisfying $\varuppsi(v^{k}a) = v^{k+1}$ and $\varuppsi(v^{i}) = v^{i}$ for $i = 1,\ldots,k$. Thus, setting $\varupphi_{k+1} := \varuppsi \circ \varupphi_{k}$ we obtain an automorphism in $\mathrm{Aut}_{0}(\T)$ satisfying $\varupphi_{k+1}(u^{i}) = v^{i}$ for $i = 1,\ldots,k{+}1$. To see that the arising sequence $(\varupphi_{k})_{k=1}^{\infty}$ has the Cauchy property note that by the properties of flip automorphisms for any $k \leq \ell$ the automorphisms $\varupphi_{k}$ and $\varupphi_{\ell}$ act identically on the ball of radius $k$, which implies $d(\varupphi_{k},\varupphi_{\ell}) \leq e^{-k}$ and shows the claim. \hfill $\Box$ 

\medskip

The final lemma shows that for any two horospheres we can find an automorphism mapping one to the other in a way that is consistent with the parametrization by $\mathbb{G}$.

\begin{theorem} \lemma \label{automorphism_horosphere}
For any two boundary points $\xi,\zeta \in \partial\T$ there is an automorphism $\varupphi \in \mathrm{Aut}_{0}(\T)$ mapping $s(g,\xi)$ to $s(g,\zeta)$ for every $g \in \mathbb{G}$.

\end{theorem}

\textit{Proof:} For $n \in \mathbb{N}$, $1 \leq j \leq n$ and $w \in \{1,\ldots,(d{-}1)\}^{j}$ with $w_{1} \neq (d{-}1)$ we define the sets $C_{n,j}^{w} \subseteq \mathbb{G}$ by
	\[C_{n,j}^{w} := g_{n}^{w_{1}}\cdots g_{n-j+1}^{w_{j}}\mathbb{G}_{n-j}\]
and set
    \[D_{n,j}^{w,\xi} := \big\{s(g,\xi)\colon g \in C_{n,j}^{w}\big\}\]
for every $\xi \in \partial\T$. Now fix two boundary points $\xi,\zeta \in \partial\T$. We shall construct an automorphism $\varupphi \in \mathrm{Aut}_{0}(\T)$ fixing $e$ and mapping the set $D_{n,j}^{w,\xi}$ to the set $D_{n,j}^{w,\zeta}$ for every $n \in \mathbb{N}$, $1 \leq j \leq n$ and $w \in \{1,\ldots,(d{-}1)\}^{j}$ with $w_{1} \neq (d{-}1)$. Note that the sets $C_{n,n}^{w}$ are singletons containing the element $g_{n}^{w_{1}}\cdots g_{1}^{w_{n}}$ and every element $g \in \mathbb{G}\setminus\{1\}$ is of this form. This implies that for every $g \in \mathbb{G}\setminus\{1\}$ there is an index $n$ and a corresponding word $w$ such that $D_{n,n}^{w,\xi}$ is the singleton containing $s(g,\xi)$ and $D_{n,n}^{w,\zeta}$ is the singleton containing $s(g,\zeta)$. Furthermore we have $s(1,\xi) = e = s(1,\zeta)$. Thus, once we have constructed the above automorphism, the statement of the theorem will follow as a special case. To prove the claim we shall define a Cauchy sequence $(\varupphi_{n})_{n=1}^{\infty}$ of automorphisms $\varupphi_{n} \in \mathrm{Aut}_{0}(\T)$ such that $\varupphi_{n}$ fixes $e$, maps the element $\xi_{1}\ldots\xi_{k}$ to the element $\zeta_{1}\ldots\zeta_{k}$ for every $k \in \mathbb{N}$ and the set $D_{m,j}^{w,\xi}$ to the set $D_{m,j}^{w,\zeta}$ for $m \leq n$, $1 \leq j \leq m$ and $w \in \{1,\ldots,(d{-}1)\}^{j}$ with $w_{1} \neq (d{-}1)$. Noting that the sets $D_{n,j}^{w,\xi}$ are contained in the ball $B_{2n}$, we may argue as in the proof of Lemma \ref{automorphism_geodesic} that the limit $\varupphi$ of the sequence $(\varupphi_{n})_{n=1}^{\infty}$ satisfies the desired property.

We shall define the sequence $(\varupphi_{n})_{n=1}^{\infty}$ inductively as follows. Let $\varupphi_{0}$ be an automorphism mapping the geodesic corresponding to $\xi$ to the geodesic corresponding to $\zeta$ as in Lemma \ref{automorphism_geodesic}. For $n = 1$ we may note that the set $C_{1,1}^{i}$ is a singleton consisting of the element $g_{1}^{i}$ for every $i \in \{1,\ldots,(d{-}2)\}$. Thus there are letters $a,b \in \Sigma$ compatible with $\xi_{1}$ and $\zeta_{1}$, respectively, such that the sets $D_{1,1}^{i,\xi}$ and $D_{1,1}^{i,\zeta}$ are singletons consisting of elements of the form $\xi_{1}a$ and $\zeta_{1}b$ for every $i \in \{1,\ldots,(d{-}2)\}$. Accordingly we may use $(d{-}2)$ flip automorphisms, whose concatenation we denote by $\varuppsi$, to ensure that $\varupphi_{1} := \varuppsi \circ \varupphi_{0}$ maps $D_{1,1}^{i,\xi}$ to $D_{1,1}^{i,\zeta}$ for $i = 1,\ldots,(d{-}2)$. Now fix $n > 1$ and assume that we have defined $\smash{\varupphi}_{n-1}$. In order to define $\varupphi_{n}$ we shall define inductively automorphisms $\varupphi_{n,j}$ for $j = 1,\ldots,n$, which agree with $\varupphi_{n-1}$ on the ball $B_{2n-2}$, such that $\varupphi_{n,j}$ maps $\smash{D_{n,k}^{w,\xi}}$ to $\smash{D_{n,k}^{w,\zeta}}$ for $k = 1,\ldots, j$ and $w \in \{1,\ldots,(d{-}1)\}^{k}$ with $w_{1} \neq (d{-}1)$ and set $\varupphi_{n} := \varupphi_{n,n}$. We start with $j = 1$ and note that the sets $D_{n,1}^{i,\zeta}$ are contained in mutually disjoint subtrees under $\zeta_{1}\ldots\zeta_{n}$ for $i = 1,\ldots,(d{-}2)$. Furthermore, the images of the sets $D_{n,1}^{i,\xi}$ under $\varupphi_{n-1}$ are contained in mutually disjoint subtrees under $\zeta_{1}\ldots\zeta_{n}$  for $i = 1,\ldots,(d{-}2)$  by the fact that $\varupphi_{n-1}$ is a tree automorphism mapping $\xi_{1}\ldots\xi_{k}$ to $\zeta_{1}\ldots\zeta_{k}$ for $k \in \{n{-}1,n,n{+}1\}$. We may thus apply $(d{-}2)$ flip automorphisms, whose concatenation we denote by $\varuppsi_{n,1}$, to rearrange these images in a way that $D_{n,1}^{i,\xi}$ is mapped to $D_{n,1}^{i,\zeta}$ for $i = 1,\ldots,(d{-}2)$ by $\varupphi_{n,1} := \varuppsi_{n,1} \circ \varupphi_{n-1}$. Now assume that we have defined $\varupphi_{n,j}$ and fix $w \in \{1,\ldots,(d{-}1)\}^{j}$ with $w_{1} \neq (d{-}1)$. Then by assumption $D_{n,j}^{w,\xi}$ is mapped to $D_{n,j}^{w,\zeta}$ by $\varupphi_{n,j}$. Accordingly, the images of the sets $D_{n,j+1}^{wi,\xi}$ are contained in mutually disjoint subtrees under $s(g_{n}^{w_{1}}\cdots g_{n-j+1}^{w_{j}},\zeta)$ for $i = 1,\ldots,(d{-}1)$. Noting that sets $D_{n,j+1}^{wi,\zeta}$ are contained in mutually disjoint subtrees under $s(g_{n}^{w_{1}}\cdots g_{n-j+1}^{w_{j}},\zeta)$  we may apply $(d{-}1)$ flip automorphisms, whose concatenation we denote by $\varuppsi_{n,j}^{w}$ to rearrange these images in a way that $D_{n,j+1}^{wi,\xi}$ is mapped to $D_{n,j+1}^{wi,\zeta}$ for $i = 1,\ldots,(d{-}1)$ by $\varuppsi_{n,j}^{w} \circ \varupphi_{n,j}$. Thus, setting 
		\[\varupphi_{n,j+1} := \prod_{w}\varuppsi_{n,j}^{w} \circ \varupphi_{n,j}\]
we obtain an automorphism mapping the set $D_{n,j+1}^{v,\xi}$ to the set $D_{n,j+1}^{v,\zeta}$ for all words $v \in \{1,\ldots,(d{-}1)\}^{j+1}$ with $v_{1} \neq (d{-}1)$. 

To see that $(\varupphi_{n})_{n=1}^{\infty}$ is a Cauchy sequence it suffices to note that for every $k \leq \ell$ the maps $\varupphi_{\ell}$ and $\varupphi_{k}$ agree on the ball of radius $k$. Noting that we have only used automorphisms in $\mathrm{Aut}_{0}(\T)$ throughout the construction it is also not difficult to show that $\varupphi_{n} \in \mathrm{Aut}_{0}(\T)$ for every $n \in \mathbb{N}$. \hfill $\Box$

\bigskip
    
	Now consider a process $(\varupalpha_{s})_{s \in S}$ taking values in a finite set $E$. We assume that $(\varupalpha_{s})_{s \in S}$ is given in canonical form and denote by $(\Omega,\mathscr{A},\varupmu)$ its configuration space. For any boundary point $\xi \in \partial\T$ we shall consider the process $(\varupbeta^{\xi}_{g})^{~}_{g \in \mathbb{G}}$ consisting of the partitions  defined by $\varupbeta^{\xi}_{g} = \{\varupbeta^{\xi,i}_{g}\}$ defined by
        \[\varupbeta^{\xi}_{g} := \varupalpha_{s(g,\xi)}\]
    for $g \in \mathbb{G}$. The next proposition shows that invariance of $(\varupalpha_{s})_{s \in S}$ under parity preserving tree automorphisms implies stationarity of the process $(\varupbeta^{\xi}_{g})^{~}_{g \in \mathbb{G}}$ for every $\xi \in \partial\T$.
    
    \begin{theorem} \proposition \label{stationarity}
        If the process $(\varupalpha_{s})_{s \in S}$ is $\mathrm{Aut}_{0}(\T)$-invariant, then the processes $(\varupbeta^{\xi}_{g})^{~}_{g \in \mathbb{G}}$ are stationary for every $\xi \in \partial\T$.
    \end{theorem}
     \textit{Proof:} Fix $\xi \in \partial\T$ and $h \in \mathbb{G}$. By Lemma \ref{automorphism_horosphere} there is an automorphism $\varupphi \in \mathrm{Aut}_{0}(\T)$ mapping $s(g,\xi)$ to $s(g,h^{-1}\xi)$ for every $g \in \mathbb{G}$. Thus, denoting by $\varuptheta$ the tree automorphism induced by $s(h,\xi)$, we obtain
	\begin{align} \label{eq2}
    s(hg,\xi) &= \varupvarphi(g^{-1}h^{-1},\xi)^{-1} = \varupvarphi(h^{-1},\xi)^{-1}\varupvarphi(g^{-1},h^{-1}\xi)^{-1} \nonumber \\ &= s(h,\xi)s(g,h^{-1}\xi) = (\varuptheta \circ \varupphi)(s(g,\xi)),
\end{align}
and therefore 
	\[\varupbeta_{hg}^{\xi} =  \varupalpha_{(\uptheta \circ \upphi)(s(g,\xi))}\]
for every $g \in \mathbb{G}$. Noting that $\varuptheta \in \mathrm{Aut}_{0}(\T)$, the definition of $(\varupbeta^{\xi}_{g})_{g \in \mathbb{G}}$ and the invariance of $(\varupalpha_{s})^{~}_{s \in S}$ with respect to $\mathrm{Aut}_{0}(\T)$ imply that the laws of the processes $(\varupbeta^{\xi}_{g})^{~}_{g \in \mathbb{G}}$ and $\smash{(\varupbeta^{\xi}_{hg})^{~}_{g \in \mathbb{G}}}$ are the same. Since $h$ was arbitrary, this shows the claim. \hfill $\Box$

\bigskip

For the next theorem we will make use of the fact that every measure preserving and ergodic action of $\mathrm{Aut}_{0}(\T)$ on a standard probability space is automatically mixing. This is a consequence of a fundamental result due to Lubotzky and Mozes, cf. \cite{LM92}, stating that the group $\mathrm{Aut}_{0}(\T)$ enjoys the Howe-Moore property, i.e.\@ every unitary representation of $\mathrm{Aut}_{0}(\T)$ without non-trivial invariant vectors has vanishing matrix coefficients. Another proof of the above property due to Pemantle can be found in \cite{Pemantle92}.

\begin{theorem} \thm \label{mixing}
Assume that the law of the process $(\varupalpha_{s})_{s \in S}$ is invariant and ergodic under the action of $\mathrm{Aut}_{0}(\T)$. Then for every $\xi \in \partial\T$ the process $(\varupbeta_{g}^{\xi})_{g \in \mathbb{G}}^{~}$ is mixing.
\end{theorem}

\textit{Proof:} By \cite[Corollary 1]{LM92} the law of $(\varupalpha_{s})_{s \in S}$ is mixing under the action of $\mathrm{Aut}_{0}(\T)$. Consider a sequence $(h_{n})_{n=1}^{\infty}$ in $\mathbb{G}$ with $h_{n} \to \infty$. We may assume without restriction that $h_{k} \neq h_{\ell}$ for $k \neq \ell$. For every $n \in \mathbb{N}$ let $\varupphi_{n}$ denote a tree automorphism in $\mathrm{Aut}_{0}(\T)$ sending $s(g,\xi)$ to $s(g,h_{n}^{-1}\xi)$ for every $g \in \mathbb{G}$, which exists by Lemma \ref{automorphism_horosphere}. Furthermore, let $\varuptheta_{n}$ denote the automorphism in $\mathrm{Aut}_{0}(\T)$ corresponding to the element $s(h_{n},\xi)$ and set $\varuppsi_{n} := \varuptheta^{~}_{n} \circ \varupphi_{n}$. We claim that $\varuppsi_{n} \to \infty$ in $\mathrm{Aut}(\T)$. To see this note that $s(h_{k},\xi) \neq s(h_{\ell},\xi)$ and thus $\varuptheta_{k}(e) \neq \varuptheta_{\ell}(e)$. Since $\varupphi_{k}(e) = e = \varupphi_{\ell}(e)$ for any $k,\ell \in \mathbb{N}$ this implies that $\varuppsi_{k}(e) \neq \varuppsi_{\ell}(e)$ and thus $d(\varuppsi_{k},\varuppsi_{\ell}) = 1$ for every $k,\ell \in \mathbb{N}$. Accordingly, the sequence $(\varuppsi_{n})_{n=1}^{\infty}$ has no converging subsequence and thus has to leave every compact set eventually. Now let $U,V \subseteq \mathbb{G}$ be finite sets and set
    \[\widehat{U}^{\xi} := \{s(g,\xi)\colon g \in U\}\]
and
    \[\widehat{V}^{\xi} := \{s(g,\xi)\colon g \in V\}.\]
Then by definition we have 
    \[(\varupbeta^{\xi}_{g})^{~}_{g \in U} = (\varupalpha_{s})_{s \in \widehat{U}^{\xi}}\]
and we may argue as in the proof of Proposition \ref{stationarity} to obtain that
    \[(\varupbeta^{\xi}_{g})_{g \in h_{n}V} = (\varupalpha_{s})_{s \in \uppsi_{n}(\widehat{V}^{\xi})}\]
for every $n \in \mathbb{N}$, see equality~\eqref{eq2}. Since $(\varupalpha_{s})_{s \in S}$ is mixing with respect to $\mathrm{Aut}_{0}(\T)$ this implies that for every $\varepsilon > 0$ there is some $n_{0} \in \mathbb{N}$ such that for all $n \geq n_{0}$ we have
    \[\sup_{i,j}\Big|\varupmu\big(\varupbeta^{\xi,i}_{h_{n}V} \cap \varupbeta^{\xi,j}_{U}\big) - \varupmu\big(\varupbeta^{\xi,i}_{V}\big)\varupmu\big(\varupbeta^{\xi,j}_{U}\big)\Big| = \sup_{i,j}\big|\varupmu\big(\varupalpha_{\uppsi_{n}(\smash{\widehat{V}^{\xi}})}^{~\!\!i} \cap \varupalpha_{U^{\xi}}^{~\!\!j}\big) - \varupmu\big(\varupalpha_{~\smash{\widehat{V}^{\xi}}}^{~\!\!i}\big)\varupmu\big(\varupalpha_{U^{\xi}}^{~\!\!j}\big)\big| < \varepsilon.\]
This shows that the process $(\varupbeta_{g}^{\xi})_{g \in \mathbb{G}}^{~}$ is mixing. \hfill $\Box$  

\bigskip

The above results allow us to apply the equipartition theory for amenable groups to the processes $(\varupbeta_{g}^{\xi})_{g \in \mathbb{G}}^{~}$. Recall that a sequence $(F_{n})_{n=1}^{\infty}$ of finite subsets of a countable group $\Gamma$ is called a \textit{F\o lner sequence} if it satisfies
    \[\lim_{n \to \infty}\frac{|F_{n} \triangle gF_{n}|}{|F_{n}|} = 0\]
for every $g \in \Gamma$ and \textit{tempered} if there is a constant $C > 0$ such that
    \[\left|\bigcup_{i=1}^{n-1}F_{i}^{-1}F_{n}\right| \leq C|F_{n}|\]
for all $n \in \mathbb{N}$. We will say that $(F_{n})_{n=1}^{\infty}$ has \textit{superlogarithmic growth} if the sequence $\log n~|F_{n}|^{-1}$ converges to zero for $n \to \infty$. 

\begin{theorem} \thm \label{ConvergenceFolner}
    Assume that the law of the process $(\varupalpha_{s})_{s \in S}$ is invariant and ergodic under the action of $\mathrm{Aut}_{0}(\T)$. Then there is a constant $h \in [0,\log|E|]$ such that the following holds. For any $\xi \in \partial\T$ let $(F^{\xi}_{n})_{n=1}^{\infty}$ be a F\o lner sequence in $\mathbb{G}$ and consider the sequence $(\widehat{F}_{n}^{\xi})_{n=1}^{\infty}$ given by
        \[\widehat{F}_{n}^{\xi} := \big\{s(g,\xi)\colon g \in F_{n}^{\xi}\big\}\]
    for $n \in \mathbb{N}$. Then we obtain
        \[\lim_{n \to \infty}\frac{I\big(\varupalpha_{\smash{\widehat{F}_{n}^{\xi}}}\big)}{|\widehat{F}_{n}^{\xi}|} = h \]
    in $L^1$ for every $\xi \in \partial\T$. If the sequences $(F^{\xi}_{n})_{n=1}^{\infty}$ are tempered and of superlogarithmic growth, we obtain as well almost sure convergence.
\end{theorem}

\textit{Proof:} By Proposition \ref{stationarity} and Theorem \ref{mixing} the processes $(\varupbeta_{g}^{\xi})_{g \in \mathbb{G}}^{~}$ are stationary and mixing, thus in particular ergodic. Fixing $\xi \in \partial\T$ we may therefore apply Kieffer's version of the Shannon-McMillan theorem for stationary and ergodic processes over amenable groups, cf. \cite{Kieffer75}, to the process $(\varupbeta_{g}^{\xi})_{g \in \mathbb{G}}^{~}$ obtain a constant $h_{\xi} \in [0,\log|E|]$ such that for any F\o lner sequence $(F_{n})_{n=1}^{\infty}$ we have
    \[\lim_{n \to \infty}\frac{I\big(\varupalpha_{\smash{\widehat{F}_{n}}}\big)}{|\widehat{F}_{n}|} = \lim_{n \to \infty}\frac{I\big(\varupbeta^{\xi}_{~\!\!F_{n}}\big)}{|F_{n}|} = h_{\xi}\]
in $L^{1}$, where we set
    \[\widehat{F}_{n} := \{s(g,\xi)\colon g \in F_{n}\}\]
for $n \in \mathbb{N}$. If the sequence $(F_{n})_{n=1}^{\infty}$ is moreover tempered and of superlogarithmic growth we obtain almost sure convergence by Lindenstrauss' generalization of the Shannon-McMillan-Breiman theorem for amenable groups, cf. \cite{Lindenstrauss01}. 

Now fix $\xi,\zeta \in \partial\T$ and let $\varupphi \in \mathrm{Aut}_{0}(\T)$ denote a tree automorphism mapping $s(g,\xi)$ to $s(g,\zeta)$ for every $g \in \mathbb{G}$ according to Lemma \ref{automorphism_horosphere}. By the invariance of $(\varupalpha_{s})_{s \in S}$ with respect to $\mathrm{Aut}_{0}(\T)$ this implies that the law of the processes $(\varupbeta_{g}^{\xi})_{g \in \mathbb{G}}^{~}$ and $(\varupbeta_{g}^{\zeta})_{g \in \mathbb{G}}^{~}$ is the same, so we have $h_{\xi} = h_{\zeta}$. This gives the claim. \hfill $\Box$

\bigskip

We shall obtain Theorem \ref{main1} as a corollary of the above theorem. To this end we consider the sequences $(U_{n})_{n=1}^{\infty}$ and $(V_{n})_{n=1}^{\infty}$ given by $U_{n} := \mathbb{G}_{n}$ and $V_{n} := U_{n} \setminus U_{n-1}$ for $n \in \mathbb{N}$. Fix $\xi \in \partial\T$. Denoting by $B_{n}^{\xi} := B_{2n} \cap H(\xi)$ and $S_{n}^{\xi} := S_{2n} \cap H(\xi)$ the horospherical balls and spheres respectively we have
    \[B_{n}^{\xi} = \{s(g,\xi)\colon g \in U_{n}\}\]
and
    \[S_{n}^{\xi} = \{s(g,\xi)\colon g \in V_{n}\}\]
for $n \in \mathbb{N}$. Therefore Theorem \ref{main1} will be proven as soon as we have shown the following lemma.

\begin{theorem} \lemma
    The sequences $(U_{n})_{n=1}^{\infty}$ and $(V_{n})_{n=1}^{\infty}$ are tempered F\o lner sequences in $\mathbb{G}$ of superlogarithmic growth.
\end{theorem}

    \textit{Proof:} For any $g \in \mathbb{G}$ we find some $n_{0} \in \mathbb{N}$ such that $g \in \mathbb{G}_{n_{0}}$. Since $\mathbb{G}_{n} \subseteq \mathbb{G}_{n+1}$ for $n \in \mathbb{N}$ this implies $g\mathbb{G}_{n} = \mathbb{G}_{n}$ for $n \geq n_{0}$. This readily implies the F\o lner property of $(U_{n})_{n=1}^{\infty}$. The temperedness of $(U_{n})_{n=1}^{\infty}$ follows from the fact that $\mathbb{G}_{i}^{-1}\mathbb{G}_{n} = \mathbb{G}_{n}$ for every $n \in \mathbb{N}$. Finally, by Proposition \ref{Proposition}, we have $|\mathbb{G}_{n}| = (d{-}1)^{n}$, so the superlogarithmic growth condition is satisfied as well.

    To show the F\o lner property for the sequence $(V)_{n=1}^{\infty}$ it suffices to note that by the same argument as above we obtain
        \[g V_{n} = g(U_{n}\setminus U_{n-1}) = gU_{n}\setminus g U_{n-1} = U_{n} \setminus U_{n-1} = V_{n}\]
    for any $g \in \mathbb{G}$ and every $n \in \mathbb{N}$ sufficiently large. Furthermore, for every $n \in \mathbb{N}$ we have
        \[|V_{n}| = (d{-}1)^{n} - (d{-}1)^{n-1} = (d{-}2)(d{-}1)^{n-1},\]
    which shows the superlogarithmic growth condition. Finally, using the above formula together with the fact that $V_{i} \subseteq U_{i}$ for every $i \in \mathbb{N}$ we obtain 
        \begin{align*}
            \left|\bigcup_{i=1}^{n-1}V_{i}^{-1}V_{n}\right| \leq \left|\bigcup_{i=1}^{n-1}U_{i}U_{n}\right| = |U_{n}| = (d{-}1)^{n} = \frac{d{-}1}{d{-}2}|V_{n}|,
        \end{align*}
    which gives the temperedness of the sequence. \hfill $\Box$
\bigskip
    
\section{Exponential mixing and equipartition along spheres}

Our aim in this section is to prove a spherical equipartition theorem for processes on regular trees satisfying a suitable uniform mixing condition, which can be defined for any process over a countable metric space. To this end, let us fix a countable infinite set $S$ and let $d$ be a metric on $S$. Consider a stochastic process $(\varupalpha_{s})_{s \in S}$ taking values in finite set $E$. Fixing two finite sets $U,V \subseteq S$ we set
    \[\varuppsi(U,V) := \sup_{i,j}\left|\frac{\varupmu(\varupalpha_{U}^{i} \cap \varupalpha_{V}^{j})}{\varupmu(\varupalpha^{i}_{U})\varupmu(\varupalpha^{j}_{V})} - 1 \right|\]
For $\lambda > 0$ we shall say that (the law of) the process $(\varupalpha_{s})_{s \in S}$ is \textit{$\varuppsi$-mixing with exponential decay rate $\lambda$} if there is a $C > 0$ such that for all finite sets $U,V \subseteq S$ we have 
    \[\varuppsi(U,V) \leq C~\!|U|~\!|V|~e^{-\lambda~\!\!d(U,V)},\]
whenever the right hand side is less or equal than $1$. Our aim is to prove the following theorem, which has been stated already as Theorem~\ref{main2} in the introduction.

\begin{theorem} \thm \label{MainTheorem}
       Let $\T = (S,B)$ be the $d$-regular tree. Consider a process $(\varupalpha_{s})_{s \in S}$ with finite state space $E$ whose law is invariant under the group $\mathrm{Aut}_{0}(\T)$ and $\varuppsi$-mixing with exponential decay rate $\lambda > 2\log(d{-}1)$. Then there exists a constant $h \in [0,\log|E|]$ such that
        \[\lim_{n \to \infty}\frac{I\big(\varupalpha_{S_{2n}}\big)}{|S_{2n}|} = h\]
    in $L^{1}$ and almost surely. 
    \end{theorem}

The above mixing condition may appear very strong a priori. While it is obviously satisfied for i.i.d.\@ processes, it may seem unclear at first if one may expect to find any interesting example beyond this setting. However, it turns out that for Gibbs fields on bounded degree graphs, exponential $\varuppsi$-mixing occurs naturally in many contexts. In the case of the $d$-regular tree it follows from results of the second author in \cite{Zimmermann25+} the unique Gibbs field realizing the Ising model at sufficiently low inverse temperature is $\varuppsi$-mixing with arbitrarily high exponential decay rate, cf. \cite[Corollary 4.3]{Zimmermann25+}. Moreover, for sufficiently large $N$ the antiferromagnetic $N$-states Potts is $\varuppsi$-mixing with arbitrarily high exponential decay rate, independently of the inverse temperature, cf. \cite[Corollary 4.4]{Zimmermann25+}. Moreover, one can show that all of the above Gibbs fields are homogeneous Markov chains, and thus, they are invariant under all tree automorphisms, see \cite[Chapter 12]{Geo11}. Consequently, they provide a rich class of examples that satisfy the conditions of theorem \ref{MainTheorem}. Note that in \cite{Zimmermann25+} a slightly stronger definition of $\varuppsi$-mixing is used, which implies the property of $\varuppsi$-mixing as defined above as a special case. 

\medskip

In the proof of Theorem \ref{MainTheorem}, we will need the following lemma. It shows that for finite state processes on regular graphs, exponential $\varuppsi$-mixing implies mixing in the ergodic theoretic sense.

\begin{theorem} \lemma   \label{psi-mixing_mixing}
    Let $G = (S,B)$ be a regular graph. Consider a process $(\varupalpha_{s})_{s \in S}$ taking values in a finite set $E$. If the law of $(\varupalpha_{s})_{s \in S}$ is exponentially $\varuppsi$-mixing, then it is also mixing under the action of $\mathrm{Aut}(G)$.
\end{theorem}

\textit{Proof:} Let $(\smash{\varupphi}_{n})_{n=1}^{\infty}$ be a sequence of automorphisms $\varupphi_{n} \in \mathrm{Aut}(T)$ with $\varupphi_{n} \to \infty$. We will show that for any sites $s,t \in S$ we have $d(s,\smash{\varupphi}_{n}t) \to \infty$ for $n \to \infty$. For the sake of a contradiction let us assume the latter statement is not true. Then there is some $k \in \mathbb{N}$ such that for infinitely many $n \in \mathbb{N}$ we have
    \[d(t,\smash{\varupphi}_{n}t) \leq d(t,s) + d(s,\smash{\varupphi}_{n}t) < k.\]
This means that $t$ is mapped into the ball $B_{k}(t)$ by infinitely many $\smash{\varupphi}_{n}t$. By the pigeonhole principle this implies that there is some $u \in S$ such that $\smash{\varupphi}_{n}t = u$ for infinitely many $n \in \mathbb{N}$. Accordingly infinitely many automorphisms $\smash{\varupphi}_{n}t$ lie in the cylinder set
    \[B := \{\varuppsi \in \mathrm{Aut}(G)\colon \varuppsi(t) = u\}.\]
It is not difficult to see that $B$ is closed (thus complete) and totally bounded. Accordingly $B$ is a compact set. This contradicts the assumption. 

Now fix $U,V \subseteq S$ finite. Then by the above observation we obtain $d(U,\smash{\varupphi}_{n}V) \to \infty$ and thus $\varuppsi(U,\varupphi_{n}V) \to 0$. Thus, for any $\varepsilon > 0$ we find some $n_{0} \in \mathbb{N}$ such that for all $n \geq n_{0}$ we have 
    \[\big|\varupmu\big(\varupalpha_{\upphi_{n}V}^{~\!\!i} \cap \varupalpha_{U}^{~\!\!j}\big) - \varupmu\big(\varupalpha_{\upphi_{n}V}^{~\!\!i}\big)\varupmu\big(\varupalpha_{U}^{~\!\!j}\big)\big| \leq \left|\frac{\varupmu\big(\varupalpha_{\upphi_{n}V}^{i} \cap \varupalpha_{U}^{j}\big)}{\varupmu\big(\varupalpha^{i}_{\upphi_{n}V}\big)\varupmu\big(\varupalpha^{j}_{U}\big)} - 1 \right| < \varepsilon\]
for all $i,j$. This shows the claim. \hfill $\Box$

\bigskip
    
In order to prove Theorem \ref{MainTheorem} we shall rely on the results of the foregoing section. We define for every $\xi \in \partial\T$ a sequence $(F_{n}^{\xi})_{n=1}^{\infty}$ of sets $F_{n}^{\xi} \subseteq \mathbb{G}$ as follows. Fixing $n \in \mathbb{N}$ and $\xi \in \partial\T$ let $a$ be the minimal letter in $\Sigma \setminus \{\xi_{n}\}$ that is larger than $\xi_{n+1}^{-1}$, if such a letter exists and the least letter in $\Sigma \setminus \{\xi_{n}\}$ otherwise. Then choose a $j \in \{1,\ldots,d{-}1\}$ such that the $n$-th letter of $g_{n}^{j}(\xi)$ equals $a$ and define
        \[F_{n}^{\xi} := g_{n}^{j}\mathbb{G}_{n-1}.\]
It is not hard to see that the set $F_{n}^{\xi}$, as well as the set $\widehat{F}_{n}^{\xi}$, are determined by the initial segment $\xi_{1}\ldots\xi_{n+1}$ of $\xi$. In the following, we will therefore write $F_{n}^{u}$ and $\widehat{F}_{n}^{u}$ instead of $F_{n}^{\xi}$ and $\widehat{F}_{n}^{\xi}$ respectively for a sequence $\xi$ starting with the initial segment $u = u_{1}\ldots u_{n+1}$. The properties of the sets $\widehat{F}_{n}^{u}$ described by the next lemma will play a crucial role in the proof of Theorem \ref{MainTheorem}.

\begin{theorem} \lemma \label{Lemma1}
    Fix $n \in \mathbb{N}$. Then the family of sets $\{\widehat{F}_{n}^{u}\}$ has the following properties.
		\begin{enumerate}[~~~(i)]
			\item For all $u,v \in S_{n+1}$ with $u \neq v$ the sets $\widehat{F}^{u}_{n}$ and $\widehat{F}^{v}_{n}$ have distance at least $2n$.
			\item For every $u \in S_{n+1}$ the set $\widehat{F}^{u}_{n}$ has $(d{-}1)^{n-1}$ elements.
			\item The family of sets $\{\widehat{F}^{u}_{n}\}$ builds a disjoint partition of the sphere $S_{2n}$.
\end{enumerate}
\end{theorem}

	\textit{Proof:} Fix $u,v \in S_{n+1}$ with $u \neq v$ and let $a$ denote the minimal letter such that $a\neq u_{n}$, which is larger than $u_{n+1}^{-1}$, if such a letter exists, and the least letter $a \neq u_{n}$ otherwise. Fix $g \in F_{n}^{u}$. By definition the element $g$ maps a sequence $\xi$ with $\xi_{1}\ldots\xi_{n+1} = u$ to a sequence $\zeta$ with $\zeta_{n} = a$. Thus, every element of $\widehat{F}^{u}_{n}$ is of the form 
	\[u_{1}\ldots u_{n} a^{-1} w_{n-1}^{-1}\ldots w_{1}^{-1}\]
for some word $w_{1}\ldots w_{n-1}$ satisfying $w_{n-1} \neq a^{-1}$. Now fix $v \in S_{n+1}$ with $v \neq u$ and let $b$ be the letter playing the same role for $v$ as $a$ does for $u$. Let $i$ be the minimal index in $\{1,\ldots,n{+}1\}$ such that $u_{i} \neq v_{i}$. Fix $s \in \widehat{F}^{u}_{n}$ and $t \in \widehat{F}^{v}_{n}$. Then $s$ is a word of length $2n$ starting in $u_{1}\ldots u_{n}$ such that $s_{n+1} = a^{-1}$ and $t$ is a word of length $2n$ starting in $v_{1}\ldots v_{n}$ such that $t_{n+1} = b^{-1}$. If $i \leq n$, we obtain
	\[d(u,v) = 2(2n{-}i{+}1) > 2n.\]
On other hand, if $i = n{+}1$, then it is not difficult to see that $a \neq b$ and thus 
    \[s_{n+1} = a^{-1} \neq b^{-1} = t_{n+1},\]
in which case we obtain $d(u,v) = 2n$. This shows statement $(i)$. Statement $(ii)$ follows from the fact that 
    \[|F_{n}^{u}| = |\mathbb{G}_{n-1}| = (d{-}1)^{n-1}\]
for $j = 1,\ldots,d{-}1$ and the map $s(\cdot,\xi)$ is injective for every $\xi \in \partial\T$. Finally, observing that $\widehat{F}^{u}_{n} \subseteq S_{2n}$ and the sets $D^{u}_{n}$ are pairwise disjoint by $(i)$, statement $(iii)$ follows from $(ii)$ by cardinality reasons. \hfill $\Box$

\begin{theorem} \lemma \label{Lemma2}
	For every $\xi \in \partial\T$ the sequence $(F^{\xi}_{n})_{n=1}^{\infty}$ is a tempered F\o lner sequence of superlogarithmic growth.
\end{theorem}
	\textit{Proof:} Fix $\xi \in \partial\T$. Then, using the fact that $\mathbb{G}$ is abelian, we obtain 
		\[|g^{j}_{n}\mathbb{G}_{n-1} \triangle hg^{j}_{n}\mathbb{G}_{n-1}| = |\mathbb{G}_{n-1} \triangle h\mathbb{G}_{n-1}|. \]
for every $j \in \{1,\ldots,d{-}1\}$, $h \in \mathbb{G}$ and $n \in \mathbb{N}$. Accordingly we get
        \[|F_{n}^{\xi} \triangle hF_{n}^{\xi}| = |\mathbb{G}_{n-1}\triangle h\mathbb{G}_{n-1}|\]
for every $n \in \mathbb{N}$ and every $h \in \mathbb{G}$. Thus the F\o lner property of $(F_{n}^{\xi})_{n=1}^{\infty}$ follows from the F\o lner property of $(\mathbb{G}_{n})_{n=1}^{\infty}$. For the temperedness fix $n \in \mathbb{N}$. Then $F^{\xi}_{n} \subseteq \mathbb{G}_{n}$ and $|\mathbb{G}_{n}| = (d{-}1) |F^{\xi}_{n}|$. This implies 
    \[\left|\bigcup_{i=0}^{n-1}\big(F^{\xi}_{i}\big)^{-1}F^{\xi}_{n}\right| \leq |\mathbb{G}_{n}| = (d{-}1)|F^{\xi}_{n}|.\]
Finally, the growth property follows from Lemma \ref{Lemma2} (ii). \hfill $\Box$ 

\medskip

The proof of the following Lemma follows the lines of Krengel's proof of Lemma 2.6 in \cite{Krengel85}. It states that for arbitrary finite state processes, the normalized information function satisfies an $L^1$-maximal inequality along sequences of finite sets with mild growth conditions.

\begin{theorem} \lemma \label{maximal_inequality}
    Let $S$ be a countable index set and $E$ be a finite set of states. Then there is a constant $C > 0$ depending only on the cardinality of $E$ such that the following holds. Consider a process $(\varupalpha_{s})_{s \in S}$ with finite state space $E$ and let $(F_{n})_{n=1}^{\infty}$ be a sequence of non-empty finite sets $F_{n} \subseteq S$ satisfying $|F_{n+1}| > |F_{n}|$ for every $n \in \mathbb{N}$. Then we have
        \[\int\sup_{n \in \mathbb{N}} \frac{I\big(\varupalpha_{F_{n}}\big)}{|F_{n}|} < C.\]
\end{theorem}

\textit{Proof:} Without restriction we may assume that the process $(\varupalpha_{s})_{s \in S}$ is given in canonical form. Let $(\Omega,\varupmu)$ be the configuration space. Fix $n \in \mathbb{N}$ and $r > 0$ and consider the set
    \[B_{r}^{n} := \{I\big(\varupalpha_{F_{n}}\big) > r|F_{n}|\}.\]
Then $B_{r}^{n}$ is the union of those atoms $B$ in the partition $\varupalpha_{F_{n}}$ satisfying
    \[\varupmu(B) < e^{-|F_{n}|r}.\]
Since the total number of atoms in $\varupalpha_{F_{n}}$ is $|E|^{|F_{n}|}$, we thus obtain 
    \[\varupmu(B_{r}^{n}) \leq |E|^{|F_{n}|}\cdot e^{-|F_{n}|r} = e^{(\log|E| - r)|F_{n}|}.\]
We now consider the set 
    \[A_{r} := \left\{\sup_{n \in \mathbb{N}}  \frac{I\big(\varupalpha_{F_{n}}\big)}{|F_{n}|}   > r\right\}.\]
It is not hard to see that
    \[A_{r} = \bigcup_{n=1}^{\infty} B_{r}^{n}.\]
Fix a number $r_{0} > \log|E|$ (depending only on $|E|$) such that for all $r > r_{0}$ we have
    \[e^{r-\log|E|} \geq 1 + e^{-2\log|E|}e^{r}.\]
Now let $r > r_{0}$. Noting that by assumption, we have $|F_{n}| \geq n$ for every $n \in \mathbb{N}$, we obtain
    \[\varupmu(A^{r}) \leq \sum_{n=1}^{\infty}e^{(\log|E| - r)n} = \frac{1}{e^{r - \log|E|} - 1} \leq e^{2\log|E|}e^{-r}.\]
By Cavalieri's principle this implies
    \begin{align*}
        \int_{\Omega}\sup_{n \in \mathbb{N}} \frac{I\big(\varupalpha_{F_{n}}\big)}{|F_{n}|}~\!\mathrm{d}\varupmu &= \int_{0}^{\infty} \varupmu(A_{r})~\!\mathrm{d} r \\ &= \int_{0}^{r_{0}}\varupmu(A_{r}) ~\!\mathrm{d} r+ \int_{r_{0}}^{\infty} \varupmu(A_{r})~\!\mathrm{d} r  \\ &\leq r_{0} + e^{2\log|E|}e^{-r_{0}}. \vphantom{\int} 
    \end{align*}
Since $r_{0}$ does only depend on $|E|$, this finishes the proof.
\hfill $\Box$

\medskip

With the above lemmas at hand, the proof of Theorem \ref{MainTheorem} proceeds as follows. By Lemma \ref{Lemma2} in combination with Lemma \ref{psi-mixing_mixing} and Theorem \ref{ConvergenceFolner}, there is a constant $h \in [0,\log|E|]$ such that for every $\xi \in \partial\T$ we have
	\[\lim_{n \to \infty} \frac{I\big(\varupalpha_{\smash{\widehat{F}_{n}^{\xi}}}\big)}{|\widehat{F}_{n}^{\xi}|} = h\]
$\varupmu\text{-almost surely}$. This implies that for $\varupmu$-almost every $\omega \in \Omega$ we obtain 
    \[\lim_{n \to \infty} \frac{I\big(\varupalpha_{\smash{\widehat{F}_{n}^{~\!\!\cdot}}}\big)(\omega)}{|\widehat{F}_{n}^{~\!\!\cdot}|} = h\]
$\varupnu$-almost surely. Now by Lemma \ref{maximal_inequality}  there is a constant $C_{1} > 0$ such that
    \[\int_{\partial T}\int_{\Omega}~\sup_{n \in \mathbb{N}}\frac{I\big(\varupalpha_{\smash{\widehat{F}_{n}^{\xi}}}\big)(\omega)}{|\widehat{F}_{n}^{\xi}|}~\mathrm{d}\varupmu(\omega)\mathrm{d}\varupnu(\xi) < C_{1}.\]
With Fubini's theorem, this implies that for $\varupmu$-almost all $\omega \in \Omega$ we have
    \[\int_{\partial\T}\sup_{n \in \mathbb{N}}\frac{I\big(\varupalpha_{\smash{\widehat{F}_{n}^{\xi}}}\big)(\omega)}{|\widehat{F}_{n}^{\xi}|}~\mathrm{d}\varupnu(\xi) < \infty.\]
Accordingly we may apply dominated convergence to obtain that for $\varupmu$-almost every $\omega \in \Omega$ we have 
	\[ \lim_{n \to \infty}~ \int_{\partial\T} \frac{I\big(\varupalpha_{\smash{\widehat{F}_{n}^{\xi}}}\big)(\omega)}{|\widehat{F}_{n}^{\xi}|}~\mathrm{d}\varupnu(\xi) = h.\]

Now fix $n \in \mathbb{N}$ and $\omega \in \Omega$ as above and set $m := d(d{-}1)^{n}$. Considering an enumeration $u^{1},\ldots,u^{m}$ of the elements of the sphere $S_{n+1}$ we set $\widehat{F}_{n}^{i} := \widehat{F}_{n}^{u^{i}}$ for $1 \leq i \leq m$ and $\widehat{F}_{n}^{i,j} := \widehat{F}_{n}^{i} \cup \ldots \cup \widehat{F}_{n}^{j}$ for $1 \leq i \leq j \leq m$. Furthermore, for $1 \leq i \leq m$, we denote by $\varupalpha_{i,n}^{\omega}$ the unique atom of $\varupalpha_{\smash{\widehat{F}_{n}^{i}}}$ containing $\omega$. Denoting furthermore by $\varupalpha_{n}^{\omega}$ the atom of $\varupalpha_{S_{2n}}$ containing $\omega$ we obtain $\varupalpha_{n}^{\omega} = \varupalpha_{1,n}^{\omega} \cap \ldots \cap \varupalpha_{m,n}^{\omega}$ by Lemma \ref{Lemma1} (iii), which gives
\begin{align*}
		I\big(\varupalpha_{S_{2n}}\big)(\omega) =  -\log\varupmu\Bigg(\bigcap_{j=1}^{m} \varupalpha^{\omega}_{j,n}\Bigg).
	\end{align*}
Furthermore, using Lemma \ref{Lemma1} (ii), we obtain 
	\begin{align*}
		\int_{\partial\T} \frac{I\big(\varupalpha_{\smash{\widehat{F}_{n}^{\xi}}}\big)(\omega)}{|\widehat{F}_{n}^{\xi}|}~\mathrm{d}\varupnu(\xi) &= \sum_{i=1}^{m} \frac{\varupnu(u_{i})}{|\widehat{F}_{n}^{i}|} ~\!I\big(\varupalpha_{\smash{\widehat{F}_{n}^{i}}}\big)(\omega) \\ &= -\frac{1}{|S_{2n}|}\sum_{i=1}^{m}\log\varupmu\big(\varupalpha_{n,i}^{\omega}\big) \\
&= -\frac{1}{|S_{2n}|}\log \prod_{i = 1}^{m} \varupmu\big(\varupalpha_{n,i}^{\omega}\big).
	\end{align*}
Accordingly, we get
\begin{align*}
	 \left| \int_{\partial\T} \frac{I\big(\varupalpha_{\smash{\widehat{F}_{n}^{\xi}}}\big)(\omega)}{|\widehat{F}_{n}^{\xi}|}~\mathrm{d}\varupnu(\xi) - \frac{I\big(\varupalpha_{S_{2n}}\big)(\omega)}{|S_{2n}|}\right|= \frac{1}{|S_{2n}|}\left|\log \frac{\varupmu\big(\bigcap_{i = 1}^{m}\varupalpha_{n,i}^{\omega}\big)}{\prod_{i = 1}^{m} \varupmu\big(\varupalpha_{n,i}^{\omega}\big)}\right|. 
\end{align*}
Now we have 
\begin{align*}
\frac{\varupmu\big(\bigcap_{i = 1}^{m}\varupalpha_{n,i}^{\omega}\big)}{\prod_{i = 1}^{m} \varupmu\big(\varupalpha_{n,i}^{\omega}\big)} &= \prod_{j=1}^{m-1} \frac{\varupmu\big(\bigcap_{i = 1}^{j+1}\varupalpha_{n,i}^{\omega}\big)}{\varupmu\big(\bigcap_{i = 1}^{j}\varupalpha_{n,i}^{\omega}\big)\varupmu\big(\varupalpha_{n,j+1}^{\omega}\big)} \leq \prod_{j=1}^{m-1} \big(1 + \varuppsi(\widehat{F}_{n}^{1,j},\widehat{F}_{n}^{j+1})\big).
\end{align*}
Let $C$ be the constant witnessing the exponential $\varuppsi$-mixing of $\varupmu$ and set $C_{0} := Cd^{2}(d{-}1)^{-2}$. Fix $\varepsilon > 0$ such that $\lambda = (2+\varepsilon)\log(d{-}1)$. Using Lemma \ref{Lemma1} $(i)$ we observe that
    \begin{align*}
    C ~|\widehat{F}_{n}^{1,j}| ~|\widehat{F}_{n}^{j+1}| e^{-\lambda d(F_{n}^{1,j},F_{n}^{j+1})}
    &\leq C|S_{2n}|^{2}e^{-(2+\varepsilon)\log(d-1) \cdot 2n} \\ &= Cd^{2}(d{-}1)^{4n-2}(d{-}1)^{-2(2 + \varepsilon) n} \\ &= C_{0} (d{-}1)^{-2\varepsilon n}.
    \end{align*}
Since for $n$ large enough the last term becomes smaller then $1$, we may apply the mixing condition to obtain that
\begin{align*}
\smash{\varuppsi(\widehat{F}_{n}^{1,j},\widehat{F}_{n}^{j+1})} &\leq C_{0} (d{-}1)^{-2\varepsilon n}
\end{align*}
for $j = 1,\ldots, m$ and thus 
\[\frac{1}{|S_{n+1}|}\log\frac{\varupmu\big(\bigcap_{i = 1}^{m}\varupalpha_{n,i}^{\omega}\big)}{\prod_{i = 1}^{m} \varupmu\big(\varupalpha_{n,i}^{\omega}\big)} \leq \log\big(1 + C_{0} (d{-}1)^{-2\varepsilon n}\big)\]
for all $n$ sufficiently large. This implies
    \[\limsup_{n \to \infty} \frac{1}{|S_{n+1}|}\log\frac{\varupmu\big(\bigcap_{i = 1}^{m}\varupalpha_{n,i}^{\omega}\big)}{\prod_{i = 1}^{m} \varupmu\big(\varupalpha_{n,i}^{\omega}\big)} \leq 0.\]
Similarly, we obtain
\begin{align*}
\frac{1}{|S_{n+1}|}\log\frac{\varupmu\big(\bigcap_{i = 1}^{m}\varupalpha_{n,i}^{\omega}\big)}{\prod_{i = 1}^{m} \varupmu\big(\varupalpha_{n,i}^{\omega}\big)} \geq \log\left(1 - C_{0} (d{-}1)^{-2\varepsilon n}\right)
\end{align*}
for all $n$ large enough and thus  
    \[\liminf_{n \to \infty} \frac{1}{|S_{n+1}|}\log\frac{\varupmu\big(\bigcap_{i = 1}^{m}\varupalpha_{n,i}^{\omega}\big)}{\prod_{i = 1}^{m} \varupmu\big(\varupalpha_{n,i}^{\omega}\big)} \geq 0.\]
Both together implies 
    \[\lim_{n \to \infty}\frac{1}{|S_{2n}|}\left|\log\frac{\varupmu\big(\bigcap_{i = 1}^{m}\varupalpha_{n,i}^{\omega}\big)}{\prod_{i = 1}^{m} \varupmu\big(\varupalpha_{n,i}^{\omega}\big)}\right| \leq \lim_{n \to \infty}\frac{1}{|S_{n+1}|}\left|\log\frac{\varupmu\big(\bigcap_{i = 1}^{m}\varupalpha_{n,i}^{\omega}\big)}{\prod_{i = 1}^{m} \varupmu\big(\varupalpha_{n,i}^{\omega}\big)}\right| = 0.\]
In sum we have shown
	\[\lim_{n \to \infty}\frac{1}{|S_{2n}|}I(\varupalpha_{S_{2n}})(\omega) = h.\]
This shows almost everywhere convergence. To verify $L^{1}$-convergence note that by Lemma \ref{maximal_inequality} the function
    \[\sup_{n \in \mathbb{N}}\left|\frac{I(\varupalpha_{S_{2n}})}{|S_{2n}|} - h\right| \leq \sup_{n \in \mathbb{N}}\frac{I(\varupalpha_{S_{2n}})}{|S_{2n}|} + h\]
is integrable with respect to $\varupmu$. Thus, we may apply dominated convergence to obtain
    \[\lim_{n \to \infty}\int_{\Omega}\left|\frac{I(\varupalpha_{S_{2n}})}{|S_{2n}|} - h\right|~\mathrm{d}\varupmu = 0.\]
This finishes the proof of Theorem \ref{MainTheorem}.
   
\end{spacing}

\textsc{Mathematical Institute, University of Leipzig}

\textsc{Augustusplatz 10, 04109 Leipzig}

\textit{felix.pogorzelski@math.uni-leipzig.de}

\textit{elias.zimmermann@math.uni-leipzig.de}

\end{document}